\documentclass[12pt,reqno]{amsart}
\usepackage{amsmath,amssymb,latexsym,esint,cite,mathrsfs}
\usepackage{verbatim,wasysym}
\usepackage[left=2.6cm,right=2.6cm,top=3cm,bottom=3cm]{geometry}
\makeatletter \@addtoreset{equation}{section} \makeatother

\usepackage{graphicx}
\usepackage{subfigure}
\usepackage{graphicx,epic,eepic}

\usepackage{color,enumitem,graphicx}
\usepackage[colorlinks=true,urlcolor=blue,
citecolor=red,linkcolor=blue,linktocpage,pdfpagelabels,
bookmarksnumbered,bookmarksopen]{hyperref}

\newtheorem{theo}{Theorem}

\numberwithin{equation}{section}
\newtheorem{theorem}{Theorem}[section]
\newtheorem{lemma}[theorem]{Lemma}
\newtheorem{definition}[theorem]{Definition}

\newtheorem{corollary}[theorem]{Corollary}
\newtheorem{remark}[theorem]{Remark}
\numberwithin{equation}{section}

\allowdisplaybreaks

\begin{document}

\title[Stability of Rellich-Sobolev inequality]
{Stability of Rellich-Sobolev type inequality involving Hardy term for bi-Laplacian}

\author[S. Deng]{Shengbing Deng}
\address{\noindent Shengbing Deng \newline
School of Mathematics and Statistics, Southwest University,
Chongqing 400715, People's Republic of China}\email{shbdeng@swu.edu.cn}

\author[X. Tian]{Xingliang Tian$^{\ast}$}
\address{\noindent Xingliang Tian  \newline
School of Mathematics and Statistics, Southwest University,
Chongqing 400715, People's Republic of China.}\email{xltian@email.swu.edu.cn}

\thanks{$^{\ast}$ Corresponding author}



\allowdisplaybreaks

\begin{abstract}
{\tiny For $N\geq 5$ and $0<\mu<N-4$, we first show a non-degenerate result of the extremal functions for the following Rellich-Sobolev type inequality
    \begin{align*}
    \int_{\mathbb{R}^N}|\Delta u|^2 \mathrm{d}x
    -C_{\mu,1}\int_{\mathbb{R}^N}\frac{|\nabla u|^2}{|x|^2} \mathrm{d}x
    +C_{\mu,2}\int_{\mathbb{R}^N}\frac{u^2}{|x|^4} \mathrm{d}x
    \geq \mathcal{S}_\mu\left(\int_{\mathbb{R}^N}|u|^{\frac{2N}{N-4}} \mathrm{d}x\right)^\frac{N-4}{N},\quad \forall u\in C^\infty_0(\mathbb{R}^N),
    \end{align*}
where $C_{\mu,1}$, $C_{\mu,2}$ and $\mathcal{S}_\mu$ are constants depending on $N$ and $\mu$, which is a key ingredient in analyzing the blow-up phenomena of solutions to various elliptic equations on bounded or unbounded domains. Then by using spectral analysis combined with a compactness argument, we consider the stability of this inequality. Furthermore, we derive a remainder term inequality in the weak Lebesgue-norm sense in a subdomain with finite Lebesgue measure.
    }
\end{abstract}

\vspace{3mm}

\maketitle

\section{{\bfseries Introduction}}\label{sectir}

    \subsection{Motivation} Let us recall the classical Sobolev inequality which states that for $N\geq 3$, there exists $\mathcal{S}=\mathcal{S}(N)>0$ such that
\begin{equation}\label{bsic}
\|\nabla u\|^2_{L^2(\mathbb{R}^N)}\geq \mathcal{S}\|u\|^2_{L^{2^*}(\mathbb{R}^N)},\quad  \mbox{for all}\ u\in \mathcal{D}^{1,2}_0(\mathbb{R}^N),
\end{equation}
where $2^*=2N/(N-2)$ and $\mathcal{D}^{1,2}_0(\mathbb{R}^N)$ denotes the completion of $C^\infty_0(\mathbb{R}^N)$ with respect to the norm $\|u\|_{\mathcal{D}^{1,2}_0(\mathbb{R}^N)}=\|\nabla u\|_{L^2(\mathbb{R}^N)}$. By using rearrangement methods, Talenti \cite{Ta76} found the optimal constant and the extremals for inequality \eqref{bsic}. Indeed, equality in \eqref{bsic} is achieved by the functions
\[V_{z,\lambda}(x)=A\left(\frac{\lambda}{1+\lambda^2|x-z|^2}
\right)^{\frac{N-2}{2}},\]
with $A\in\mathbb{R}$, $z\in\mathbb{R}^N$ and $\lambda>0$. It is well known that the Euler-Lagrange equation associated to (\ref{bsic}) is
\begin{equation}\label{bec}
-\Delta u=|u|^{2^*-2}u,\quad \mbox{in}\ \mathbb{R}^N.
\end{equation}
By Caffarelli et al. \cite{CGS89}, it is known that all positive solutions of \eqref{bec} are Talenti bubble $V_{z,\lambda}(x)$ with $A=[N(N-2)]^{\frac{N-2}{4}}$.
In \cite{BrL85}, Brezis and Lieb asked the question whether a remainder term - proportional to the quadratic distance of the function $u$ to be the manifold  $\mathcal{M}:=\{cV_{\lambda,z}: c\in\mathbb{R}, \lambda>0, z\in\mathbb{R}^N\}$ - can be added to the right hand side of \eqref{bsic}. This question was answered affirmatively by Bianchi and Egnell \cite{BE91}, which reads that there exists constant $c_{\mathrm{BE}}>0$ such that for all $u\in \mathcal{D}^{1,2}_0(\mathbb{R}^N)$,
    \begin{equation}\label{defcbe}
    \|\nabla u\|^2_{L^2(\mathbb{R}^N)}
    -\mathcal{S}\|u\|^2_{L^{2^*}(\mathbb{R}^N)}
    \geq c_{\mathrm{BE}} \inf\limits_{v\in \mathcal{M}}\|\nabla u-\nabla v\|^2_{L^2(\mathbb{R}^N)}.
    \end{equation}
    Furthermore, for the weighted case, R\u{a}dulescu et al. \cite{RSW02} gave the remainder term of Hardy-Sobolev inequality for exponent two (see also our recent work \cite{DTp} involving $p$-Laplace), and recently, Wei and Wu \cite{WW22} established the stability of the profile decompositions to the Caffarelli-Kohn-Nirenberg inequality and also gave the remainder term.

    Then, let us recall the well-known second-order Sobolev inequality (see \cite{EFJ90,Li85-1,Va93})
    \begin{equation}\label{bcesi}
    \|\Delta u\|^2_{L^2(\mathbb{R}^N)}
    \geq \mathcal{S}_0\|u\|^2_{L^{2^{**}}(\mathbb{R}^N)},\quad  \forall u\in \mathcal{D}^{2,2}_0(\mathbb{R}^N),
    \end{equation}
    with the sharp constant
    \begin{equation}\label{cssi}
    \mathcal{S}_0=\pi^2N(N-4)(N^2-4)
    \left(\frac{\Gamma(N/2)}{\Gamma(N)}\right)^{4/N},
    \end{equation}
    where $N\geq 5$, $2^{**}:=\frac{2N}{N-4}$, $\mathcal{D}^{2,2}_0(\mathbb{R}^N)$ denotes the completion of $C^\infty_0(\mathbb{R}^N)$ with respect to the norm $\|u\|_{\mathcal{D}^{2,2}_0(\mathbb{R}^N)}=\|\Delta u\|_{L^2(\mathbb{R}^N)}$, and $\Gamma$ denotes the Gamma function $\Gamma(\gamma):=\int^{\infty}_0 t^{\gamma-1} e^{-t}\mathrm{d}t$ for $\gamma>0$. Note that the related Euler-Lagrange equation is
    \begin{equation}\label{bce}
    \Delta^2 u=|u|^{2^{**}-2}u,\quad \mbox{in}\ \mathbb{R}^N.
    \end{equation}
    Smooth positive solutions of (\ref{bce}) have been completely classified by  Lin \cite{Li98}, that is, the author proved that they are given by $\lambda^{\frac{N-4}{2}}U_0(\lambda(x-z))$ for $\lambda>0$ and $z\in\mathbb{R}^N$, where
    \begin{equation*}
    U_0(x)=[(N-4)(N-2)N(N+2)]^{\frac{N-4}{8}}
    \left(1+|x|^2\right)^{-\frac{N-4}{2}},
    \end{equation*}
    and they are unique (up to scalar multiplications) optimizers for inequality \eqref{bcesi}.
    In \cite{LW00}, Lu and Wei considered the following linearized equation
    \begin{equation*}
    \Delta^2 v=\xi U_0^{2^{**}-2}v,\quad \mbox{in}\ \mathbb{R}^N,\quad v\in \mathcal{D}^{2,2}_0(\mathbb{R}^N),
    \end{equation*}
    and they proved that the eigenvalues $\{\xi\}$ of the above problem are discrete:
    \[
    \xi_1=1,\quad \xi_2=\xi_3=\cdots=\xi_{N+2}=2^{**}-1<\xi_{N+3}\leq\cdots,
    \]
    and the corresponding eigenfunction spaces are
    \[
    V_1=\mathrm{Span}\{U_0\},\quad V_2=\mathrm{Span}\left\{\frac{N-4}{2}U_0+x\cdot\nabla U_0, \quad  \frac{\partial U_0}{\partial x_i},i=1,\ldots,N\right\}, \cdots.
    \]
    Moreover, Lu and Wei \cite{LW00} proved the stability of inequality \eqref{bcesi} which extends the work of Bianchi and Egnell \cite{BE91} to the second-order case, i.e., there is $c_{\mathrm{LW}}>0$ such that
    \begin{equation}\label{bcesir}
    \|\Delta u\|^2_{L^2(\mathbb{R}^N)}
    -\mathcal{S}_0\|u\|^2_{L^{2^{**}}(\mathbb{R}^N)}
    \geq c_{\mathrm{LW}}\inf\limits_{v\in \mathcal{M}_0}\|\Delta(u-v)\|^2_{L^2(\mathbb{R}^N)},\quad \forall u\in \mathcal{D}^{2,2}_0(\mathbb{R}^N),
    \end{equation}
    where $\mathcal{M}_0:=\{c\lambda^{\frac{N-4}{2}}U_0(\lambda(x-z)): c\in\mathbb{R},\ \lambda>0,\ z\in\mathbb{R}^N\}$ is the set of extremal functions for inequality \eqref{bcesi}. We also refer to \cite{BWW03,CFW13,FZ22} for the general higher order cases about the stability of Sobolev inequality.

    Let us recall the Rellich inequality (see \cite{Re69}): for $N\geq 5$,
    \begin{align}\label{Ri}
    \int_{\mathbb{R}^N}|\Delta u|^2 \mathrm{d}x
    \geq \left(\frac{N(N-4)}{4}\right)^2\int_{\mathbb{R}^N}\frac{u^2}{|x|^4} \mathrm{d}x,\quad \forall u\in C^\infty_0(\mathbb{R}^N),
    \end{align}
    the constant $\left(\frac{N(N-4)}{4}\right)^2$ is sharp, and the inequality is strict for any nontrivial functions.
    Then it is easy to deduce the Rellich-Sobolev inequality: for $N\geq 5$, $0\leq \alpha<\left(\frac{N(N-4)}{4}\right)^2$ and $0\leq \beta<4$,
    \begin{align*}
    \mathcal{S}_{\alpha,\beta}
    :=\inf_{u\in \mathcal{D}^{2,2}_0(\mathbb{R}^N)\setminus\{0\}}
    \frac{\int_{\mathbb{R}^N}|\Delta u|^2 \mathrm{d}x
    -\alpha\int_{\mathbb{R}^N}\frac{u^2}{|x|^4} \mathrm{d}x}
    {\left(\int_{\mathbb{R}^N}\frac{|u|^{2^{**}_\beta}}{|x|^\beta} \mathrm{d}x\right)^\frac{2}{2^{**}_\beta}}>0,
    \end{align*}
    where $2^{**}_\beta=\frac{2(N-\beta)}{N-4}$, and from \cite{BM12}
    we know that $\mathcal{S}_{\alpha,\beta}$ can be achieved. Furthermore, by the classical rearrangement methods, we know that the minimizers are unique (up to scalings and multiplications) for any $\alpha+\beta>0$. However, the explicit form of minimizers is not known yet and only for its asymptotic behavior is known, see \cite{DJ15,JL14,KX17} for details. Tertikas and Zographopoulos \cite{TZ07} established the following Hardy-Rellich inequality: for $N\geq 5$,
    \begin{align}\label{HRi}
    \int_{\mathbb{R}^N}|\Delta u|^2 \mathrm{d}x
    \geq \frac{N^2}{4}\int_{\mathbb{R}^N}\frac{|\nabla u|^2}{|x|^2} \mathrm{d}x,\quad \forall u\in C^\infty_0(\mathbb{R}^N),
    \end{align}
    and the constant $\frac{N^2}{4}$ is sharp, and the inequality is strict for any nontrivial functions, see also \cite{Caz20} for more general case. Then based on inequality \eqref{HRi}, Caldiroli \cite{Cal14} considered the Hardy-Rellich-Sobolev inequality: for $0< \lambda<\frac{N^2}{4}$ and $0\leq \beta<4$,
    \begin{align*}
    \mathcal{S}_{\lambda,\beta}
    :=\inf_{u\in \mathcal{D}^{2,2}_0(\mathbb{R}^N)\setminus\{0\}}
    \frac{\int_{\mathbb{R}^N}|\Delta u|^2 \mathrm{d}x
    -\lambda\int_{\mathbb{R}^N}\frac{|\nabla u|^2}{|x|^2} \mathrm{d}x}
    {\left(\int_{\mathbb{R}^N}\frac{|u|^{2^{**}_\beta}}{|x|^\beta} \mathrm{d}x\right)^\frac{2}{2^{**}_\beta}}>0,
    \end{align*}
    and the constant $\mathcal{S}_{\lambda,\beta}$ can be achieved, but the explicit form of minimizers is also not known yet, see also \cite{CC16}.

    As mentioned above, it is natural to consider whether the explicit minimizers can be given for some weighted second-order inequalities. Fortunately, Dan, Ma and Yang in \cite[Theorem 1.6]{DMY20} established the following Rellich-Sobolev type inequality with the explicit form of minimizers:

    \begin{theo}\label{thmDMY} Let $N\geq 5$ and $0<\mu<N-4$. There holds, for all $u\in \mathcal{D}^{2,2}_0(\mathbb{R}^N)$,
    \begin{align}\label{RSi}
    & \int_{\mathbb{R}^N}|\Delta u|^2 \mathrm{d}x
    -C_{\mu,1}\int_{\mathbb{R}^N}\frac{|\nabla u|^2}{|x|^2} \mathrm{d}x
    +C_{\mu,2}\int_{\mathbb{R}^N}\frac{u^2}{|x|^4} \mathrm{d}x
    \geq \mathcal{S}_\mu\left(\int_{\mathbb{R}^N}|u|^{2^{**}} \mathrm{d}x\right)^\frac{2}{2^{**}},
    \end{align}
    where
    \begin{align}\label{RSisc}
    C_{\mu,1}:=& \frac{N^2-4N+8}{2(N-4)^2}\mu[2(N-4)-\mu];
    \nonumber\\
    C_{\mu,2}:=& \frac{N^2}{16(N-4)^2}\mu^2[2(N-4)-\mu]^2
    -\frac{N-2}{2}\mu[2(N-4)-\mu];
    \nonumber\\
    \mathcal{S}_\mu:=& \left(1-\frac{\mu}{N-4}\right)^{4-\frac{4}{N}}\mathcal{S}_0,
    \end{align}
    and $\mathcal{S}_0$ is the sharp constant of second-order Sobolev inequality given in \eqref{cssi}.
    Moreover, the constant $\mathcal{S}_\mu$ in \eqref{RSi} is sharp and equality holds if and only if
    \begin{align*}
    u(x)=c|x|^{-\frac{\mu}{2}}
    \left(\lambda^2+|x|^{2(1-\frac{\mu}{N-4})}
    \right)^{-\frac{N-4}{2}},
    \end{align*}
    for $c\in\mathbb{R}$ and $\lambda>0$.
    \end{theo}

Note that when $\mu=0$, the Rellich-Sobolev type inequality \eqref{RSi} reduces to the classical second-order Sobolev inequality \eqref{bcesi}. We also refer to \cite{BM12,DGT23-jde,DT23,DT24,HW20,Ya21} for other weighted second-order inequalities which are restricted in radial space, and in those references minimizers are classified.

\subsection{Main results}
    In the present paper, we are mainly concerned about the stability of inequality \eqref{RSi}, that is, we will give its remainder term in an appropriate space which extends the work of Lu and Wei \cite{LW00} to the Hardy-Rellich case.

    Define
    \begin{align*}
    \mathcal{E}_\mu:=
    \left\{
    u\in \mathcal{D}^{2,2}_0(\mathbb{R}^N): \int_{\mathbb{R}^N}|\Delta u|^2 \mathrm{d}x
    -C_{\mu,1}\int_{\mathbb{R}^N}\frac{|\nabla u|^2}{|x|^2} \mathrm{d}x
    +C_{\mu,2}\int_{\mathbb{R}^N}\frac{u^2}{|x|^4} \mathrm{d}x
    <\infty
    \right\},
    \end{align*}
    with respect to the inner product
    \begin{align*}
    \langle u, v\rangle_\mu
    := \int_{\mathbb{R}^N}\Delta u \Delta v\mathrm{d}x
    -C_{\mu,1}\int_{\mathbb{R}^N}\frac{\nabla u\cdot\nabla v}{|x|^2} \mathrm{d}x
    +C_{\mu,2}\int_{\mathbb{R}^N}\frac{uv}{|x|^4} \mathrm{d}x,
    \end{align*}
    and the norm $\|u\|_\mu:=\langle u, u\rangle_\mu^{1/2}$.
    From Rellich inequality \eqref{Ri}, Hardy-Rellich inequality \eqref{HRi}, and \eqref{RSi}, it is not difficult to verify that the space $\mathcal{E}_\mu$ is well defined and equivalent to $\mathcal{D}^{2,2}_0(\mathbb{R}^N)$.

    We first concern the uniqueness of positive radial solutions of the Euler-Lagrange equation
    \begin{align}\label{RSe}
    & \Delta^2 u
    +C_{\mu,1}\mathrm{div}\left(\frac{\nabla u}{|x|^2}\right)
    +C_{\mu,2}\frac{u}{|x|^4}
    =|u|^{2^{**}-2}u,\quad \mbox{in}\ \mathbb{R}^N\setminus\{\mathbf{0}\},\quad u\in \mathcal{E}_\mu,
    \end{align}
    where $N\geq 5$ and $0<\mu<N-4$, $C_{\mu,1}$ and $C_{\mu,2}$ are constants depending only on $N$ and $\mu$ defined in \eqref{RSisc}.

    \begin{theorem}\label{thmprtu}
    Suppose $N\geq 5$ and $0<\mu<N-4$. Then problem \eqref{RSe} admits a unique (up to scalings) positive radial solution of the form $U_{\mu,\lambda}(x)=\lambda^{\frac{N-4}{2}}U_\mu(\lambda x)$ for $\lambda>0$, where
    \begin{equation}\label{defulaeb}
    U_\mu(x)=K_{N,\mu}|x|^{-\frac{\mu}{2}}
    \left(1+|x|^{2(1-\frac{\mu}{N-4})}
    \right)^{-\frac{N-4}{2}},
    \end{equation}
    with
    \[
    K_{N,\mu}=\left[\left(1-\frac{\mu}{N-4}
    \right)^4(N-4)(N-2)N(N+2)\right]^{\frac{N-4}{8}}.
    \]
    \end{theorem}

    Then we concern the linearized problem related to Euler-Lagrange equation \eqref{RSe} at the function $U_\mu$, which leads to consider the following problem:
    \begin{align}\label{RSel}
    & \Delta^2 v
    +C_{\mu,1}\mathrm{div}\left(\frac{\nabla v}{|x|^2}\right)
    +C_{\mu,2}\frac{v}{|x|^4}
    =(2^{**}-1)U_\mu^{2^{**}-2}v,\quad \mbox{in} \ \mathbb{R}^N\setminus\{0\}, \quad v\in \mathcal{E}_\mu.
    \end{align}
    It is obvious that $\frac{N-4}{2}U_\mu+x\cdot \nabla U_\mu$ (which equals to $\frac{\partial U_{\mu,\lambda}}{\partial \lambda}|_{\lambda=1}$) solves the linear equation \eqref{RSel}. We say that $U_\mu$ is non-degenerate if all the solutions of \eqref{RSel} result from the invariance (up to scalings) of \eqref{RSe}. The non-degeneracy of solutions for \eqref{RSe} is a key ingredient in analyzing the blow-up phenomena of solutions whose asymptotic behavior is encoded in \eqref{defulaeb}. Therefore, it is quite natural to ask the following question:
    \begin{center}
    \em Is the solution $U_{\mu}$ non-degenerate?
    \end{center}
    We give an affirmative answer in the following.

    \begin{theorem}\label{thmlpnd}
    Suppose $N\geq 5$ and $0<\mu<N-4$. Then the space of solutions for (\ref{RSel}) has dimension $1$ and is spanned by $(\frac{N-4}{2}U_\mu+x\cdot \nabla U_\mu)$.
    \end{theorem}

    A direct application of Theorem \ref{thmlpnd} is to study the stability of Rellich-Sobolev type inequality \eqref{RSi}.
    Now, let us state our second main result.

    \begin{theorem}\label{thmprt}
    Suppose $N\geq 5$ and $0<\mu<N-4$. There exists a constant  $\mathcal{C}_1=\mathcal{C}(N,\mu)>0$ such that for all $u\in \mathcal{E}_\mu$,
    \[
    \|u\|_\mu^2
    -\mathcal{S}_\mu\|u\|^2_{L^{2^{**}}(\mathbb{R}^N)}
    \geq \mathcal{C}_1\inf\limits_{v\in \mathcal{M}_\mu}\|u-v\|^2_\mu,
    \]
    where $\mathcal{M}_\mu:=\{cU_{\mu,\lambda}: c\in\mathbb{R}, \lambda>0\}$ is the set of extremal functions for Rellich-Sobolev type inequality \eqref{RSi}.
    \end{theorem}

    \begin{remark}\label{remc}\rm
    The key step of the proofs for Theorems \ref{thmprtu} and \ref{thmlpnd} is the change of variable $u(r)=r^av(r^b)$ for radial function with
    \[
    a=-\frac{\mu}{2}\quad \mbox{and}\quad b=1-\frac{\mu}{N-4},
    \]
    which requires very careful calculation. The proof of Theorem \ref{thmprt} is standard, but surprisingly, a very useful tool established by Dan et al. \cite{DMY20} is good for us, that is, under the variable $u(x)=|x|^av(|x|^{b-1}x)$ for general function,
    \begin{align*}
    \|u\|_\mu^2\geq \left(1-\frac{\mu}{N-4}\right)^3
    \|\Delta v\|^2_{L^2(\mathbb{R}^N)},
    \end{align*}
    and the equality holds if and only if $u$ is radially symmetric. Furthermore,
    \begin{align*}
    \|u\|^{2^{**}}_{L^{2^{**}}(\mathbb{R}^N)}
    =\left(1-\frac{\mu}{N-4}\right)^{-1}
    \|v\|^{2^{**}}_{L^{2^{**}}(\mathbb{R}^N)},
    \end{align*}
    thus
    \begin{align*}
    \|\Delta v\|^2_{L^2(\mathbb{R}^N)}
    -\mathcal{S}_0\|v\|^2_{L^{2^{**}}(\mathbb{R}^N)}
    \leq
    \left(1-\frac{\mu}{N-4}\right)^{-3}\left(\|u\|_\mu^2
    -\mathcal{S}_\mu\|u\|^2_{L^{2^{**}}(\mathbb{R}^N)}
    \right),
    \end{align*}
    then we can directly use Lions' concentration-compactness principle to complete the proof.
    \end{remark}

    For a domain $\Omega\subset\mathbb{R}^N$ with smooth boundary containing origin in its interior, we define
    \[
    \mathcal{S}_\mu(\Omega):=\inf_{u\in \dot{\mathcal{E}}_{\mu,\Omega}\setminus\{0\}}
    \frac{\|u\|_{\mu,\Omega}^2}
    {\|u\|_{L^{2^{**}}(\Omega)}^2}.
    \]
    Here $\dot{\mathcal{E}}_{\mu,\Omega}:=\{u\in \mathcal{E}_{\mu,\Omega}: \tilde{u}\in \mathcal{E}_{\mu}\}$, where
    $\tilde{u}$ denotes the trivial extension of a function $u\in \mathcal{E}_{\mu,\Omega}$ on $\mathbb{R}^N$ and $\mathcal{E}_{\mu,\Omega}$ is $\mathcal{E}_{\mu}$ replacing $\mathbb{R}^N$ by $\Omega$, and
    \begin{align*}
    \|u\|^2_{\mu,\Omega}=\int_{\Omega}|\Delta u|^2 \mathrm{d}x
    -C_{\mu,1}\int_{\Omega}\frac{|\nabla u|^2}{|x|^2} \mathrm{d}x
    +C_{\mu,2}\int_{\Omega}\frac{u^2}{|x|^4} \mathrm{d}x.
    \end{align*}
    Since $\Omega$ contains origin in its interior, by using blow-up analysis, it is not difficult to verify that $\mathcal{S}_\mu(\Omega)=\mathcal{S}_\mu(\mathbb{R}^N)
    =\mathcal{S}_\mu$, and $\mathcal{S}_\mu(\Omega)$ is achieved if and only if $\Omega=\mathbb{R}^N$.
    Therefore, based on the stability result of extremal functions as in Theorem \ref{thmprt}, we will derive a weak-Lebesgue remainder term inequality as \cite{BE91,BrL85,CFW13,GG01,GW10,Ge03,PR14,RSW02,WaWi03} in a subdomain  $\Omega\subset\mathbb{R}^N$ with $|\Omega|<\infty$. Here $|\Omega|$ denotes the Lebesgue measure.

    \begin{theorem}\label{thmrtbd}
    Suppose $N\geq 5$, $0<\mu<N-4$, and let $\Omega\subset \mathbb{R}^N$ with $|\Omega|<\infty$ which has smooth boundary and contains origin in its interior. Then there exists a constant $\mathcal{C}_2=\mathcal{C}(N,\mu)>0$ independent of $\Omega$ such that
    \begin{align}\label{thmprtbd}
    \|u\|_{\mu,\Omega}^2
    -\mathcal{S}_\mu
    \|u\|_{L^{2^{**}}(\Omega)}^2
    \geq \mathcal{C}_2|\Omega|^{-\frac{N-4-\mu}{N}}
    \|u\|_{L^{\bar{p}}_w(\Omega)}^2,\quad \forall u\in \dot{\mathcal{E}}_{\mu,\Omega},
    \end{align}
    where $\bar{p}=\frac{2N}{2(N-4)-\mu}$, and $\|\cdot\|_{L^{s}_w(\Omega)}$ denotes the weak $L^{s}$-norm as
    \begin{align}\label{defwln}
    \|u\|_{L^{s}_w(\Omega)}:=
    \sup\limits_{D\subset\Omega, |D|>0}|D|^{-\frac{s-1}{s}}
    \int_{D}|u|\mathrm{d}x.
    \end{align}
    \end{theorem}

    Note that the weak $L^s$-norm as in \eqref{defwln} is equivalent to the classical weak $L^s$-norm for $s>1$, i.e., $u\in L^s_w(\Omega)$ if and only if $\sup_{t>0} t|x
    \in\Omega: |u(x)|>t|^{1/s}<\infty$. Furthermore, for any $0<t<s$ with $u\in L^s_w(\Omega)$, it holds that
    \[
    \|u\|^t_{L^t(\Omega)}\leq \frac{s}{s-t}|\Omega|^{1-\frac{t}{s}}\|u\|^t_{L^s_w(\Omega)},
    \]
    see \cite[Chapter 5]{CR16} for details. Thus, as a direct corollary of Theorem \ref{thmrtbd}, we obtain a Br\'{e}zis-Nirenberg \cite{BN83} type remainder term inequality in the following.

\begin{corollary}\label{coroprtp}
     Suppose $N\geq 5$, $0<\mu<N-4$, and let $\Omega\subset \mathbb{R}^N$ with $|\Omega|<\infty$ which has smooth boundary and contains origin in its interior. Then for each $t\in (0, \frac{2N}{2(N-4)-\mu})$, there exists a constant $\mathcal{C}_3=\mathcal{C}(N,\mu,t)>0$ independent of $\Omega$ such that
    \begin{align}\label{prtp}
    \|u\|_{\mu,\Omega}^2
    -\mathcal{S}_\mu
    \|u\|_{L^{2^{**}}(\Omega)}^2
    \geq \mathcal{C}_3|\Omega|^{\frac{2}{2^{**}}-\frac{2}{t}}
    \|u\|_{L^{t}(\Omega)}^2,\quad \forall u\in \dot{\mathcal{E}}_{\mu,\Omega}.
    \end{align}
    \end{corollary}

\subsection{Structure of the paper}\label{subsect:structrue}
    In Section \ref{appsue}, we show that equation \eqref{RSe} has a unique (up to scalings) positive radial solution of the form \eqref{defulaeb}, and give the proof of Theorem \ref{thmprtu}. Section \ref{sectndr} is devoted to proving the non-degeneracy of $U_\mu$. In Section \ref{sect:rt}, we study the stability of Rellich-Sobolev inequality \eqref{RSi} by using spectral analysis combining with a compactness argument, and give the proof of Theorem \ref{thmprt}. Finally, based on the stability result Theorem \ref{thmprt}, in Section \ref{sectwlrt} we deduce a remainder term inequality in the weak Lebesgue-norm sense as in \eqref{thmprtbd}, and give the proof of Theorem \ref{thmrtbd}.

\section{Uniqueness of radial solutions}\label{appsue}

    In this section, we will show that problem \eqref{RSe} only admits a unique (up to scalings) positive radial solution of the form \eqref{defulaeb}.

    \vskip0.25cm

\noindent{\bf \em Proof of Theorem \ref{thmprtu}.}
    Let $u\in \mathcal{E}_\mu$ be a positive radial solution of \eqref{RSe}, and set $r=|x|$, then \eqref{RSe} is equivalent to
    \begin{align}\label{RSer}
    & u^{(4)}+\frac{2(N-1)}{r}u'''
    +\frac{(N-1)(N-3)+C_{\mu,1}}{r^2}u''
    \nonumber
    \\
    &\quad -\frac{(N-3)(N-1-C_{\mu,1})}{r^3}u'
    +\frac{C_{\mu,2}}{r^4}u
    =u^{\frac{N+4}{N-4}},
    \end{align}
    where $C_{\mu,1}$ and $C_{\mu,2}$ are constants depending only on $N$ and $\mu$ given as in \eqref{RSisc}.
    Making the change of variable
    \begin{align}\label{defvc}
    u(r)=r^av(r^b),
    \end{align}
    with
    \[
    a=-\frac{\mu}{2}\quad \mbox{and}\quad b=1-\frac{\mu}{N-4},
    \]
    by direct calculations, we obtain
    \begin{align*}
    u'(r) & =ar^{a-1}v+br^{a+b-1}v', \\
    u''(r) & =a(a-1)r^{a-2}v+b(2a+b-1)r^{a+b-2}v'+b^2r^{a+2b-2}v'', \\
    u'''(r) & =a(a-1)(a-2)r^{a-3}v
    +b\left[a(a-1)+(2a+b-1)(a+b-2)\right]r^{a+b-3}v'
    \\ \nonumber & \quad+b^2\left[(2a+b-1)+(a+2b-2)\right]r^{a+2b-3}v''
    +b^3r^{a+3b-3}v''', \\
    u^{(4)}(r) & =a(a-1)(a-2)(a-3)r^{a-4}v
    \\ \nonumber &\quad +b\left\{\left[a(a-1)+(2a+b-1)(a+b-2)\right](a+b-3)
    +a(a-1)(a-2)\right\}r^{a+b-4}v'
    \\ \nonumber &\quad +b^2\left\{\left[a(a-1)+(2a+b-1)(a+b-2)\right]
    +(3a+3b-3)(a+2b-3)\right\}r^{a+2b-4}v''
    \\ \nonumber &\quad +b^3\left[(3a+3b-3)+(a+3b-3)\right]r^{a+3b-4}v'''
    +b^4r^{a+4b-4}v^{(4)}.
    \end{align*}
    Then from \eqref{RSer}, set $s=r^b$, we deduce that
    \begin{align}\label{reb}
    & v^{(4)}(s)+\frac{A}{s}v'''(s)
    +\frac{B}{s^2}v''(s)
    -\frac{C}{s^3}v'(s)
    +\frac{D}{s^4}v(s)
    =b^{-4}v^{\frac{N+4}{N-4}},
    \end{align}
    where
    \begin{align*}
    A:& = b^{-1} \left\{\left[(3a+3b-3)+(a+3b-3)\right]+2(N-1)\right\}, \\
    B:& = b^{-2} \big(\left\{\left[a(a-1)+(2a+b-1)(a+b-2)\right]
    +(3a+3b-3)(a+2b-3)\right\}
    \nonumber \\
    & \quad +2(N-1)(3a+3b-3)
    + [(N-1)(N-3)+C_{\mu,1}]
    \big), \\
    C:& = b^{-3} \big(\left\{\left[a(a-1)+(2a+b-1)(a+b-2)\right](a+b-3)
    +a(a-1)(a-2)\right\}
    \nonumber \\ & \quad +2(N-1)\left[a(a-1)+(2a+b-1)(a+b-2)\right]
    \nonumber \\ & \quad +[(N-1)(N-3)+C_{\mu,1}](2a+b-1)
    -(N-3)(N-1-C_{\mu,1})
    \big), \\
    D:& = b^{-4} \big(a(a-1)(a-2)(a-3)+2(N-1)a(a-1)(a-2)
    \\ \nonumber & \quad +[(N-1)(N-3)+C_{\mu,1}]a(a-1)
    -(N-3)(N-1-C_{\mu,1})a
    +C_{\mu,2}
    \big).
    \end{align*}
    By a direct and careful calculation, we have
    \begin{align}
    A=2(N-1),\quad B=C=(N-1)(N-3),\quad D=0,
    \end{align}
    then \eqref{reb} reduces into
    \begin{align}\label{RSeba}
    & v^{(4)}(s)+\frac{2(N-1)}{s}v'''(s)
    +\frac{(N-1)(N-3)}{s^2}v''(s)
    -\frac{(N-1)(N-3)}{s^3}v'(s)
    =b^{-4}v^{\frac{N+4}{N-4}},
    \end{align}
    and also
    \begin{eqnarray*}
    \left\{ \arraycolsep=1.5pt
       \begin{array}{ll}
        -v''+\frac{N-1}{s}v'=w, \quad \mbox{in}\ s\in(0,\infty),\\[2mm]
        -w''+\frac{N-1}{s}w'=b^{-4}v^{\frac{N+4}{N-4}}, \quad \mbox{in}\ s\in(0,\infty).
        \end{array}
    \right.
    \end{eqnarray*}
    It is not difficult to verify that equation \eqref{RSeba} admits a solution of the form
    \begin{align*}
    v(s)=\frac{K_{N,\mu}\lambda^{\frac{N-4}{2}}}
    {(1+\lambda^{2}s^{2})^{\frac{N-4}{2}}},\quad \mbox{with}\ K_{N,\mu}=\left[\left(1-\frac{\mu}{N-4}
    \right)^4(N-4)(N-2)N(N+2)\right]^{\frac{N-4}{8}},
    \end{align*}
    for any $\lambda>0$. Therefore, \eqref{RSer} admits a solution of the form $u(|x|)=\lambda^{\frac{N-4}{2}}U_\mu(\lambda x)$, where
    \begin{equation*}
    U_\mu(x)=K_{N,\mu}|x|^{-\frac{\mu}{2}}
    \left(1+|x|^{2(1-\frac{\mu}{N-4})}
    \right)^{-\frac{N-4}{2}}.
    \end{equation*}

    Next, we will show the uniqueness (up to scalings) of solutions for  equation \eqref{RSeba}.
    Following the work of Huang and Wang \cite{HW20}, we use the classical Emden-Fowler transformation, i.e., let
    \begin{align}\label{l-et}
    v(s)=b^{\frac{N-4}{2}}s^{-\frac{N-4}{2}}\varphi(t),\quad w(s)=b^{\frac{N-4}{2}}s^{-\frac{N}{2}}\phi(t),\quad \mbox{with}\ t=-\ln s.
    \end{align}
    A direct computation shows that a positive radial function $u\in \mathcal{E}_\mu$ solves \eqref{RSer} if and only if the functions $\varphi, \phi$ satisfy
    \begin{eqnarray*}
    \left\{ \arraycolsep=1.5pt
       \begin{array}{ll}
        -\varphi''+2\varphi'+\frac{N(N-4)}{4}\varphi=\phi,\quad \mbox{in}\ \mathbb{R},\\[2mm]
        -\phi''-2\phi'+\frac{N(N-4)}{4}\phi
        =\varphi^{\frac{N+4}{N-4}},\quad \mbox{in}\ \mathbb{R},
        \end{array}
    \right.
    \end{eqnarray*}
    that is, $\varphi(t)$ satisfies the following fourth-order ordinary differential equation
    \begin{equation}\label{Pwht}
    \varphi^{(4)}-K_2\varphi''+K_0\varphi
    =\varphi^{\frac{N+4}{N-4}}, \quad \varphi>0,\quad \mbox{in}\  \mathbb{R},\quad \varphi\in H^2(\mathbb{R}),
    \end{equation}
    with the constants
    \begin{align*}
    K_2=\frac{(N-2)^2+4}{2}, \quad
    K_0=\frac{N^2(N-4)^2}{16}.
    \end{align*}
    Here $H^2(\mathbb{R})$ denotes the completion of $C^\infty_0(\mathbb{R})$ with respect to the norm
    \[
    \|\varphi\|_{H^2(\mathbb{R})}
    =\left(\int_{\mathbb{R}}(|\varphi_{tt}|^2+K_2|\varphi_t|^2
    +K_0 \varphi^2)\mathrm{d}t\right)^{1/2}.
    \]
    Clearly, $K_2^2\geq 4 K_0$, then applying
    \cite[Theorem 2.2]{BM12} directly, we get the uniqueness (up to translations and inversion $t\mapsto -t$) of solutions to problem \eqref{Pwht}. Therefore, \eqref{RSer} admits a unique (up to scalings) solution $\lambda^{\frac{N-4}{2}}U_\mu(\lambda x)$ for $\lambda>0$. The proof of Theorem \ref{thmprtu} is thereby completed.
    \qed

\vskip0.25cm

\section{{\bfseries Non-degenerate result}}\label{sectndr}

    In this section, we will prove the non-degeneracy of $U_\mu$, that is, we classify all solutions of linearized problem \eqref{RSel}.

    Firstly, let us consider the eigenvalue problem
    \begin{align}\label{Pwhlep}
    \Delta^2 v
    +C_{\mu,1}\mathrm{div}\left(\frac{\nabla v}{|x|^2}\right)
    +C_{\mu,2}\frac{v}{|x|^4}
    =\nu U_{\mu}^{2^{**}-2}v,\quad \mbox{in} \ \mathbb{R}^N\setminus\{\mathbf{0}\}, \quad v\in \mathcal{E}_\mu.
    \end{align}
    In order to do this, firstly we need to establish the compactness of embedding $\mathcal{E}_\mu\hookrightarrow L^2(\mathbb{R}^N, U^{2^{**}-2}_{\mu}\mathrm{d}x)$. The following Poincar\'{e} type inequalities will be used later.
    \begin{lemma}\label{lemcztj}
    Suppose $N\geq 5$ and $0<\mu<N-4$. It holds that for all $u\in \mathcal{E}_\mu$,
    \begin{align}\label{continouseb}
    \int_{\mathbb{R}^N}
    U_\mu^{2^{**}-2}u^2\mathrm{d}x
    \leq \|u\|^2_\mu,
    \end{align}
    and the equality holds if and only if $u=\zeta U_\mu$ with $\zeta\in\mathbb{R}$.
    Also, there exist $C>0$ and $\vartheta=\vartheta(N,\mu)>0$ such that, for any $\rho\in (0,1)$,
    \begin{align}\label{continouseb1}
    \int_{B_\rho}
    U_\mu^{2^{**}-2}u^2\mathrm{d}x
    \leq C\rho^{\vartheta}\|u\|^2_\mu,
    \end{align}
    and
    \begin{align}\label{continouseb2}
    \int_{\mathbb{R}^N\setminus B_{\frac{1}{\rho}}}
    U_\mu^{2^{**}-2}u^2\mathrm{d}x
    \leq C\rho^{\vartheta}\|u\|^2_\mu.
    \end{align}
    \end{lemma}

    \begin{proof}
    By using H\"{o}lder inequality and Rellich-Sobolev type inequality \eqref{RSi} we deduce that for all $u\in \mathcal{E}_\mu$,
    \begin{align*}
    \int_{\mathbb{R}^N}U_{\mu}^{2^{**}-2}u^2 \mathrm{d}x
    \leq \|U_{\mu}\|_{L^{2^{**}}(\mathbb{R}^N)}^{2^{**}-2}
    \|u\|_{L^{2^{**}}(\mathbb{R}^N)}^{2}
    = \mathcal{S}_\mu \|u\|_{L^{2^{**}}(\mathbb{R}^N)}^{2}
    \leq \|u\|^2_\mu,
    \end{align*}
    and the first equality holds if and only if $u=\zeta U_\mu$ with $\zeta\in\mathbb{R}$, thanks to $\|U_{\mu}\|^2_\mu=\|U_{\mu}\|_{L^{2^{**}}(\mathbb{R}^N)}^{2^{**}}$ and $\|U_{\mu}\|^2_\mu=\mathcal{S}_\mu \|U_{\mu}\|_{L^{2^{**}}(\mathbb{R}^N)}^{\frac{2}{2^{**}}}
    $ which imply $\|U_{\mu}\|_{L^{2^{**}}(\mathbb{R}^N)}^{2^{**}-2}=\mathcal{S}_\mu$, then \eqref{continouseb} holds.

    To prove \eqref{continouseb1}-\eqref{continouseb2}, we apply the Sobolev inequality with radial weights. By using the change of variable $u(x)=|x|^av(|x|^{b-1}x)$ as in Remark \ref{remc}, then
    \begin{align*}
    \|u\|_\mu^2
    \geq & \left(1-\frac{\mu}{N-4}\right)^3
    \int_{\mathbb{R}^N}|\Delta v|^2 \mathrm{d}x,
    \\
    \int_{B_\rho}
    U_\mu^{2^{**}-2}u^2\mathrm{d}x
    = & \left(1-\frac{\mu}{N-4}\right)^3\int_{B_{\rho^b}}
    U_0^{2^{**}-2}v^2\mathrm{d}x,
    \\
    \int_{\mathbb{R}^N\setminus B_{\frac{1}{\rho}}}
    U_\mu^{2^{**}-2}u^2\mathrm{d}x
    = & \left(1-\frac{\mu}{N-4}\right)^3\int_{\mathbb{R}^N\setminus B_{\frac{1}{\rho^b}}}
    U_0^{2^{**}-2}v^2\mathrm{d}x.
    \end{align*}
    Note that for $\rho\in (0,1)$ and $b>0$, the Rellich inequality indicates
    \begin{align*}
    \int_{B_{\rho^b}}
    U_0^{2^{**}-2}v^2\mathrm{d}x
    \leq C \rho^{4b}\int_{B_{\rho^b}}\frac{v^2}{|x|^4}\mathrm{d}x
    \leq  C \rho^{4b}\int_{\mathbb{R}^N}
    |\Delta v|^2\mathrm{d}x,
    \end{align*}
    thus \eqref{continouseb1} holds. Moreover,
    \begin{align*}
    \int_{\mathbb{R}^N\setminus B_{\frac{1}{\rho^b}}}
    U_0^{2^{**}-2}v^2\mathrm{d}x
    \leq C\int_{\mathbb{R}^N\setminus B_{\frac{1}{\rho^b}}}\frac{v^2}{|x|^8}\mathrm{d}x
    \leq C \rho^{4b}\int_{\mathbb{R}^N\setminus B_{\frac{1}{\rho^b}}}\frac{v^2}{|x|^4}\mathrm{d}x
    \leq  C \rho^{4b}\int_{\mathbb{R}^N}
    |\Delta v|^2\mathrm{d}x,
    \end{align*}
    as our desired estimate \eqref{continouseb2}. The proof is completed.
    \end{proof}

    \begin{theorem}\label{propcet}
    Suppose $N\geq 5$ and $0<\mu<N-4$. The space $\mathcal{E}_\mu$ embeds compactly into $L^2(\mathbb{R}^N, U^{2^{**}-2}_{\mu}\mathrm{d}x)$.
    \end{theorem}

    \begin{proof}
    Since $U_\mu$ is singular at origin, we follow the arguments as those in the proof of Proposition 3.2 in \cite{FZ22}.
    Let $\{u_n\}$ be a sequence of functions in $\mathcal{E}_\mu$ with uniformly bounded norm. It follows from \eqref{continouseb} that $\|u_n\|_{L^2(\mathbb{R}^N, U^{2^{**}-2}_{\mu}\mathrm{d}x)}$ are uniformly bounded as well.

    Since $U_\mu^{2^{**}-2}$ is locally bounded away from zero and infinity in $\mathbb{R}^N\setminus\{\mathbf{0}\}$, by Rellich-Kondrachov Theorem and a diagonal argument, we deduce that, up to a subsequence, there exists $u\in \mathcal{E}_\mu$ such that $u_n\rightharpoonup u$ weakly in $\mathcal{E}_\mu$ and $u_n\to u$ strongly in $L^2_{\mathrm{loc}}(\mathbb{R}^N\setminus\{\mathbf{0}\}, U^{2^{**}-2}_{\mu}\mathrm{d}x)$.
    Also, it follows from \eqref{continouseb1} and \eqref{continouseb2} that, for any $\rho\in (0,1)$,
    \begin{align*}
    \int_{B_\rho}
    U_\mu^{2^{**}-2}u_n^2\mathrm{d}x\leq  C\rho^{\vartheta},
    \quad
    \int_{\mathbb{R}^N\setminus B_{\frac{1}{\rho}}}
    U_\mu^{2^{**}-2}u_n^2\mathrm{d}x
    \leq C\rho^{\vartheta},
    \end{align*}
    for some $\vartheta>0$. We conclude the proof by defining the compact set $K_\rho:=\overline{B_{\frac{1}{\rho}}}\setminus B_\rho$ and applying the strong convergence of $u_n$ in $K_\rho$, together with the arbitrariness of $\rho$ (that can be chosen arbitrarily small).
    \end{proof}

    Theorem \ref{propcet} indicates that the eigenvalues of problem (\ref{Pwhlep}) are discrete, then following the work of Servadei and Valdinoci \cite{SV13} for nonlocal case, we can give the definitions of eigenvalues of problem (\ref{Pwhlep}) in the following.

    \begin{definition}\label{defevp}
    The first eigenvalue of problem (\ref{Pwhlep}) can be defined as
    \begin{equation}\label{deffev1}
    \nu_1:=\inf_{v\in \mathcal{E}_\mu \setminus\{0\}}
    \frac{\|v\|^2_\mu}
    {\int_{\mathbb{R}^N}U_{\mu}^{2^{**}-2}v^2 \mathrm{d}x}.
    \end{equation}
    Moreover, for any $k\in\mathbb{N}^+$ the eigenvalues can be characterized as follows:
    \begin{equation}\label{deffevk}
    \nu_{k+1}:=\inf_{v\in \mathbb{P}_{k+1}\setminus\{0\}}
    \frac{\|v\|^2_\mu}
    {\int_{\mathbb{R}^N}U_{\mu}^{2^{**}-2}v^2 \mathrm{d}x}
    ,
    \end{equation}
    where
    \begin{equation*}
    \mathbb{P}_{k+1}:=\left\{v\in \mathcal{E}_\mu: \langle v, e_{i,j}\rangle_\mu=0,\quad \mbox{for all}\quad i=1,\ldots,k,\ j=1,\ldots,h_i\right\},
    \end{equation*}
    and $e_{i,j}$ are the corresponding eigenfunctions to $\nu_i$ with $h_i$ multiplicity.
    \end{definition}

    \begin{theorem}\label{propep}
    Suppose $N\geq 5$ and $0<\mu<N-4$. Let $\nu_i$, $i=1,2,\ldots,$ be the eigenvalues of (\ref{Pwhlep}) in increasing order as in Definition \ref{defevp}. Then $\nu_1=1$ is simple with eigenfunction $U_{\mu}$ and $\nu_2=2^{**}-1$ with the corresponding one-dimensional eigenfunction space spanned by $(\frac{N-4}{2}U_\mu+x\cdot\nabla U_\mu)$. Furthermore, $\nu_3>\nu_2$.
    \end{theorem}

    \begin{proof}
    We follow the arguments as those in \cite{BWW03}. Choosing a test function $v=U_{\mu}$ in (\ref{deffev1}), since $U_{\mu}$ is the solution of equation (\ref{RSe}) we have
    \begin{equation*}
    \begin{split}
    \nu_1
    \leq \frac{\int_{\mathbb{R}^N}|\Delta U_{\mu}|^2 \mathrm{d}x
    -C_{\mu,1}\int_{\mathbb{R}^N}\frac{|\nabla U_{\mu}|^2}{|x|^2} \mathrm{d}x
    +C_{\mu,2}\int_{\mathbb{R}^N}\frac{U_{\mu}^2}{|x|^4} \mathrm{d}x }{\int_{\mathbb{R}^N}U_{\mu}^{2^{**}}\mathrm{d}x}
    =1.
    \end{split}
    \end{equation*}
    Moreover, as in \eqref{continouseb}, for all $v\in \mathcal{E}_\mu\setminus\{0\}$,
    \begin{align}\label{fev1}
    \int_{\mathbb{R}^N}U_{\mu}^{2^{**}-2}v^2 \mathrm{d}x
    \leq \|U_{\mu}\|_{L^{2^{**}}(\mathbb{R}^N)}^{2^{**}-2}
    \|v\|_{L^{2^{**}}(\mathbb{R}^N)}^{2}
    = \mathcal{S}_\mu \|v\|_{L^{2^{**}}(\mathbb{R}^N)}^{2}
    \leq \|v\|^2_\mu,
    \end{align}
    which implies $\nu_1\geq 1$, then we have $\nu_1=1$. Furthermore, since the first equality in \eqref{fev1} holds if and only if $v=\zeta U_{\mu}$ with $\zeta\in\mathbb{R}$, corresponding eigenfunction of $\nu_1=1$ is $\zeta U_{\mu}$ with $\zeta\in\mathbb{R}\setminus\{0\}$.

    Note that $U_{\mu}$ minimizes the functional
    \begin{align}\label{deffe}
    v\mapsto \Phi(v)=\frac{1}{2}\|v\|^2_\mu-\frac{1}{2^{**}}
    \|v\|^{2^{**}}_{L^{2^{**}}(\mathbb{R}^N)},
    \end{align}
    on the Nehari manifold
    \begin{align*}
    \mathcal{N}:=\left\{v\in \mathcal{E}_\mu \backslash\{0\}: \|v\|^2_\mu=\|v\|^{2^{**}}_{L^{2^{**}}(\mathbb{R}^N)}\right\}.
    \end{align*}
    Indeed, for $v\in \mathcal{N}$, we have by Rellich-Sobolev type inequality \eqref{RSi} that
    \begin{align*}
    \Phi(v)
    & =\left(\frac{1}{2}-\frac{1}{2^{**}}\right)
    \|v\|^{2^{**}}_{L^{2^{**}}(\mathbb{R}^N)}
    = \left(\frac{1}{2}-\frac{1}{2^{**}}\right)\left(\frac{\|v\|_\mu}
    {\|v\|_{L^{2^{**}}(\mathbb{R}^N)}}\right)^{\frac{2\cdot 2^{**}}{2^{**}-2}}
    \\
    & \geq\left(\frac{1}{2}-\frac{1}{2^{**}}\right)
    \mathcal{S}_\mu^{\frac{2^{**}}{2^{**}-2}}
    =  \left(\frac{1}{2}-\frac{1}{2^{**}}\right)\left(\frac{\|U_{\mu}\|_\mu}
    {\|U_{\mu}\|_{L^{2^{**}}(\mathbb{R}^N)}}\right)^{\frac{2\cdot 2^{**}}{2^{**}-2}}
    = \Phi(U_{\mu}).
    \end{align*}
    As a consequence, the second derivative $\Phi''(U_{\mu})$ given by
    \begin{align*}
    (\phi,\varphi)\mapsto \langle \phi,\varphi\rangle_\mu
    -(2^{**}-1)\int_{\mathbb{R}^N}U_{\mu}^{2^{**}-2}\phi\varphi \mathrm{d}x
    \end{align*}
    is nonnegative quadratic form when restricted to the tangent space $T_{U_{\mu}}\mathcal{N}$, then we have
    \[
    \|u\|^2_\mu
    \geq (2^{**}-1)\int_{\mathbb{R}^N}U_{\mu}^{2^{**}-2}u^2\mathrm{d}x,
    \quad \mbox{for all}\ u\in T_{U_{\mu}}\mathcal{N}.
    \]
    Since $T_{U_{\mu}}\mathcal{N}$ has codimension one, we infer that $\nu_2\geq 2^{**}-1$. Moreover, since $(\frac{N-4}{2}U_\mu+x\cdot\nabla U_\mu)$ is a solution of \eqref{Pwhlep} with $\nu=2^{**}-1$, which indicates that $\nu_2\leq 2^{**}-1$, then we conclude that $\nu_2= 2^{**}-1$. Next, we will show that $(\frac{N-4}{2}U_\mu+x\cdot\nabla U_\mu)$ is the only (up to multiplications) eigenfunction of $\nu_2= 2^{**}-1$.

    Eigenvalue problem \eqref{Pwhlep} with $\nu=2^{**}-1$ is equivalent to
    \begin{equation}\label{Pwhlp}
    \Delta^2 v
    +C_{\mu,1}\left[|x|^{-2}\Delta v-2|x|^{-4}x\cdot\nabla v\right]
    +\frac{C_{\mu,2}}{|x|^4}v
    =\frac{(2^{**}-1)K_{N,\mu}^{2^{**}-2}|x|^{-\frac{4\mu}{N-4}} }{\left(1+|x|^{2(1-\frac{\mu}{N-4})}
    \right)^4}v,
    \end{equation}
    in $\mathbb{R}^N\setminus\{\mathbf{0}\}$, $v\in \mathcal{E}_\mu$. Firstly, we write $v(x)=v(r,\theta)$ and make the standard spherical decomposition of $v$ as follows:
    \begin{equation}\label{defvd}
    v(r,\theta)=\sum^{\infty}_{k=0}\sum^{m_k}_{i=1}\phi_{k,i}(r)
    \Psi_{k,i}(\theta),
    \end{equation}
    where $r=|x|$, $\theta=\frac{x}{|x|}\in \mathbb{S}^{N-1}$, and
    \begin{equation*}
    \phi_{k,i}(r)=\int_{\mathbb{S}^{N-1}}v(r,\theta)\Psi_{k,i}(\theta)
    \mathrm{d}\theta.
    \end{equation*}
    Here $\Psi_{k,i}(\theta)$ denotes the $k$-th spherical harmonic, i.e., it satisfies
    \begin{equation*}
    -\Delta_{\mathbb{S}^{N-1}}\Psi_{k,i}=\lambda_k \Psi_{k,i},
    \end{equation*}
    where $\Delta_{\mathbb{S}^{N-1}}$ is the Laplace-Beltrami operator on $\mathbb{S}^{N-1}$ with the standard metric and  $\lambda_k$ is the $k$-th eigenvalue of $-\Delta_{\mathbb{S}^{N-1}}$. It is well known that $\lambda_k=k(N-2+k)$, $k=0,1,2,\ldots$ whose multiplicity is
    \[
    m_k:=\frac{(N+2k-2)(N+k-3)!}{(N-2)!k!}
    \]
    and that
    \[
    \mathrm{Ker}(\Delta_{\mathbb{S}^{N-1}}+\lambda_k)
    =\mathbb{Y}_k(\mathbb{R}^N)|_{\mathbb{S}^{N-1}},
    \]
    where $\mathbb{Y}_k(\mathbb{R}^N)$ is the space of all homogeneous harmonic polynomials of degree $k$ in $\mathbb{R}^N$.
    It is known that
    \begin{align*}
    \Delta (\phi_{k,i}(r)\Psi_{k,i}(\theta))
    & =\Psi_{k,i}\left(\phi''_{k,i}+\frac{N-1}{r}\phi'_{k,i}\right)
    +\frac{\phi_{k,i}}{r^2}\Delta_{\mathbb{S}^{N-1}}\Psi_{k,i} \nonumber\\
    & =\Psi_{k,i}\left(\phi''_{k,i}+\frac{N-1}{r}\phi'_{k,i}
    -\frac{\lambda_k}{r^2}\phi_{k,i}\right),
    \end{align*}
    and
    \begin{align*}
    x\cdot\nabla (\phi_{k,i}(r)\Psi_{k,i}(\theta))
    & =\sum^{N}_{i=1}x_i\frac{\partial (\phi_{k,i}(r)\Psi_{k,i}(\theta))}{\partial x_i}
    \nonumber\\
    & =\phi'_{k,i}r\Psi_{k,i}+\phi_{k,i}
    \frac{\partial\Psi_{k,i}}{\partial \theta_l}\sum^{N}_{i=1}\frac{\partial\theta_l}{\partial x_i}x_i
    =\phi'_{k,i}r\Psi_{k,i},
    \end{align*}
    because it holds true that
    \begin{align*}
    \sum^N_{i=1}\frac{\partial\theta_l}{\partial x_i}x_i=0,\quad \mbox{for all}\ l=1,\ldots,N-1.
    \end{align*}
    Therefore, by standard regularity theory, the function $v$ is a solution of (\ref{Pwhlp}) if and only if for all $k\in\mathbb{N}$, $i=1,\ldots,m_k$, $\phi_{k,i}\in \mathcal{D}_k$ is a classical solution of
    \begin{small}
    \begin{align}\label{p2c}
    &\phi_{k,i}^{(4)}+\frac{2(N-1)}{r}\phi_{k,i}'''
    +\frac{(N-1)(N-3)+C_{\mu,1}}{r^2}\phi_{k,i}''
    -\frac{(N-3)(N-1-C_{\mu,1})}{r^3}\phi_{k,i}'
    +\frac{C_{\mu,2}}{r^4}\phi_{k,i}
    \nonumber\\[3mm]
    &\quad -\frac{\lambda_k}{r^2}
        \left[
        2\phi_{k,i}''+\frac{2(N-3)}{r}\phi_{k,i}'-\frac{2(N-4)
        +\lambda_k-C_{\mu,1}}{r^2}\phi_{k,i}
        \right]
        =\frac{(2^{**}-1)K_{N,\mu}^{2^{**}-2}r^{-\frac{4\mu}{N-4}} }{\left[1+r^{2(1-\frac{\mu}{N-4})}
        \right]^4}\phi_{k,i},
    \end{align}
    \end{small}
    in $r\in(0,\infty)$, where
    \begin{align*}
    \mathcal{D}_k:=\left\{\omega\in C_0^\infty((0,\infty)): \int^\infty_0 \left[\left(\Delta_r \omega-\frac{\lambda_k}{r^2}\omega\right)^2
    -C_{\mu,1}\frac{|\omega'|^2}{r^{2}}
    +C_{\mu,2}\frac{\omega^2}{r^{4}}\right]r^{N-1}\mathrm{d}r
    <\infty\right\}.
    \end{align*}
    Here $\Delta_r=\frac{\partial^2}{\partial r^2}+\frac{N-1}{r}\frac{\partial}{\partial r}$.
    The same as in Section \ref{appsue}, we make the change $s=r^{b}$ where $b=1-\frac{\mu}{N-4}$ and let
    \begin{equation}\label{p2txy}
    \phi_{k,i}(r)=r^{-\frac{\mu}{2}}X_{k,i}(s),
    \end{equation}
    which transforms (\ref{p2c}) into the following form
    \begin{align}\label{p2ty}
    & X_{k,i}^{(4)}+\frac{2(N-1)}{s}X_{k,i}'''
    +\frac{(N-1)(N-3)}{s^2}X_{k,i}''
    -\frac{(N-1)(N-3)}{s^3}X_{k,i}'
    \nonumber\\[3mm]
    & \quad -\frac{\lambda_kb^{-2}}{s^2}
        \left[
        2X_{k,i}''+\frac{2(N-3)}{s}X_{k,i}'-\frac{2(N-4)
        +\lambda_kb^{-2}}{s^2}X_{k,i}
        \right]
    \nonumber \\
        & =\frac{(N+4)(N-2)N(N+2)}{(1+s^2)^4}X_{k,i}
        +\frac{4b^{-2}(b^{-2}-1)\lambda_k}{s^4}X_{k,i}, &
    \end{align}
    in $s\in(0,\infty)$, and
    \[
    X_{k,i}\in \widetilde{\mathcal{D}}_k:=\left\{\omega\in C_0^\infty((0,\infty)): \int^\infty_0 \left(\Delta_s \omega-
    \frac{\lambda_k}{s^2}\omega\right)^2s^{N-1}\mathrm{d}s
    <\infty\right\}.
    \]
    Here we have used the fact
    \begin{equation*}
    b^{-4}(2^{**}-1)K_{N,\mu}^{2^{**}-2}
    =\left[(N-4)(N-2)N(N+2)\right]\left(\frac{2N}{N-4}-1\right)
    =(N+4)(N-2)N(N+2).
    \end{equation*}
    It is not difficult to verify that \eqref{p2ty} is equivalent to the system
    \begin{align}\label{rwevpb}
    \left(\Delta_s-\frac{\lambda_k}{s^2}\right)^2X_{k,i}
    & = -\frac{(b^{-2}-1)\lambda_k[2(N-4)+(1+b^2)\lambda_k-4b^{-2}]}{s^4}
    X_{k,i}
    \nonumber \\
    & \quad + \frac{2(b^{-2}-1)\lambda_k}{s^2}
    \left(X''_{k,i}+\frac{N-3}{s}X'_{k,i}\right)
    \nonumber \\
    & \quad+(2^{**}-1)\Gamma_N(1+s^2)^{-4}X_{k,i},
    \end{align}
    in $s\in(0,\infty)$, $X_{k,i}\in \widetilde{\mathcal{D}}_k$, where
    \begin{align}\label{defgn}
    \Gamma_N:=(N-4)(N-2)N(N+2).
    \end{align}

    From \cite[Lemma 2.4]{BWW03}, we know that when $k=0$, \eqref{rwevpb} only admits a unique solution $X_0=(\frac{N-4}{2}U_0+x\cdot\nabla U_0)$ (up to multiplications), that is, \eqref{p2c} with $k=0$ only admits a unique solution $(\frac{N-4}{2}U_{\mu}+x\cdot\nabla U_{\mu})$ (up to multiplications). Let us explain the reason here. In fact, the conclusion of \cite[Lemma 2.4]{BWW03} states that if $X\in \widetilde{\mathcal{D}}_0$ is a solution of
    \[
    \left[s^{1-N}\frac{\partial}{\partial s}\left(s^{N-1}\frac{\partial}{\partial s}\right)\right]^2X=\nu(1+s^2)^{-4}X \quad \mbox{for}\ \nu>0,
    \]
    with $X(0)=0$, then $X\equiv 0$. Note that $s^{1-N}\frac{\partial}{\partial s}\left(s^{N-1}\frac{\partial}{\partial s}\right)=\Delta_s$. Let $\tilde{X}_0\not\equiv 0$ be another solution of \eqref{rwevpb} with $k=0$. Then $X_0(0), \tilde{X}_0(0)\neq 0$, for otherwise $X_0$ resp. $\tilde{X}_0$ would vanish identically by \cite[Lemma 2.4]{BWW03}. So we can find $\tau\in\mathbb{R}$ such that $X_0(0)=\tau\tilde{X}_0(0)$. But then $X_0-\tau\tilde{X}_0$ also solves \eqref{rwevpb} with $k=0$ and equals zero at origin. By \cite[Lemma 2.4]{BWW03}, one has $X_0-\tau\tilde{X}_0\equiv 0$, and so $\tilde{X}_0$ is a scalar multiple of $X_0$.

    {\bf We claim that for all $k\geq 1$, \eqref{rwevpb} does not exist nontrivial solutions.}
    Now, we begin to show this claim. By direct computations, one easily checks the operator identity
    \[
    \Delta_s-\frac{\lambda_k}{s^2}
    =s^k\left(\frac{\partial^2}{\partial s^2}+\frac{N+2k-1}{s}\frac{\partial}{\partial s}\right)s^{-k}.
    \]
    Therefore, equation \eqref{rwevpb} can be rewritten as
    \begin{align}\label{rwevpbb}
    & \left(\frac{\partial^2}{\partial s^2}+\frac{N+2k-1}{s}\frac{\partial}{\partial s}\right)^2 Y_{k,i}
    \nonumber \\
    & = (b^{-2}-1)\lambda_k[2k(N-4+k)-2(N-4)-(1+b^{-2})\lambda_k+4b^{-2}]
    \frac{Y_{k,i}}{s^4}
    \nonumber \\
    & \quad+ 2(b^{-2}-1)\lambda_k
    \left[\frac{Y''_{k,i}}{s^2}
    +(N-3+2k)\frac{Y'_{k,i}}{s^3}\right]
    +(2^{**}-1)\Gamma_N(1+s^2)^{-4}Y_{k,i}.
    \end{align}
    Here we defined $Y_{k,i}\in C^\infty ((0,\infty))$ by $Y_{k,i}(s):=s^{-k}X_{k,i}$.
    Now we consider the functions $Z_{k,i}: \mathbb{R}^{N+2k}\to \mathbb{R}$ defined by $Z_{k,i}(y)=Y_{k,i}(|y|)$. So following the work of Bartsch et al. \cite{BWW03}, we deduce that
    \[
    Z_{k,i}\in \mathcal{D}^{2,2}_0(\mathbb{R}^{N+2k}),
    \]
    and $Z_{k,i}$ is a weak solution of the equation
    \begin{align}\label{rwevpbbb}
    \Delta^2 Z_{k,i}(y) & = (b^{-2}-1)\lambda_k[2k(N-4+k)-2(N-4)-(1+b^{-2})\lambda_k+4b^{-2}]
    \frac{Z_{k,i}(y)}{|y|^4}
    \nonumber \\
    & \quad + 2(b^{-2}-1)\lambda_k
    \left[\frac{Z''_{k,i}}{|y|^2}
    +(N-3+2k)\frac{Z'_{k,i}}{|y|^3}\right]
    \nonumber \\
    & \quad + (2^{**}-1)\Gamma_N(1+|y|^2)^{-4}Z_{k,i}(y),\quad y\in \mathbb{R}^{N+2k}.
    \end{align}
    Multiplying \eqref{rwevpbbb} by $Z_{k,i}$ and integrating in $\mathbb{R}^{N+2k}$, we have
    \begin{align*}
    & \|Z_{k,i}\|^2_{\mathcal{D}^{2,2}_0(\mathbb{R}^{N+2k})}
    \\
    & =(2^{**}-1)\Gamma_N\int_{\mathbb{R}^{N+2k}}
    (1+|y|^2)^{-4}|Z_{k,i}|^2 \mathrm{d}y
    -2(b^{-2}-1)\lambda_k
    \int_{\mathbb{R}^{N+2k}}
    \frac{|\nabla Z_{k,i}|^2}{|y|^{2}}
    \mathrm{d}y
    \nonumber\\
    & \quad+ (b^{-2}-1)\lambda_k
    [2k(N-4+k)-2(N-4)-(1+b^{-2})\lambda_k+4b^{-2}]
    \int_{\mathbb{R}^{N+2k}}\frac{|Z_{k,i}|^2}
    {|y|^{4}}\mathrm{d}y.
    \end{align*}
    By the Hardy inequality,
    \begin{align*}
    \left(\frac{N+2k-4}{2}\right)^2\int_{\mathbb{R}^{N+2k}}\frac{|u|^2}
    {|y|^{4}}\mathrm{d}y
    \leq \int_{\mathbb{R}^{N+2k}}
    \frac{|\nabla u|^2}{|y|^{2}}
    \mathrm{d}y, \quad \mbox{for all}\ u\in C^\infty_0(\mathbb{R}^{N+2k}),
    \end{align*}
    we deduce that
    \begin{align}\label{rwevpbbbcb}
    \|Z_{k,i}\|^2_{\mathcal{D}^{2,2}_0(\mathbb{R}^{N+2k})}
    & \leq - \frac{1}{2}(b^{-2}-1)\lambda_k
    [N(N-4)+(2\lambda_k-8)b^{-2}+2\lambda_k]
    \int_{\mathbb{R}^{N+2k}}\frac{|Z_{k,i}|^2}
    {|y|^{4}}\mathrm{d}y
    \nonumber\\
    & \quad +(2^{**}-1)\Gamma_N\int_{\mathbb{R}^{N+2k}}
    (1+|y|^2)^{-4}|Z_{k,i}|^2 \mathrm{d}y.
    \end{align}
    Note that
    \begin{align}\label{rwevpbbbcbi}
    \|u\|^2_{\mathcal{D}^{2,2}_0(\mathbb{R}^{N+2k})}
    \geq \Gamma_{N+2k}\int_{\mathbb{R}^{N+2k}}
    (1+|y|^2)^{-4}|u(y)|^2 \mathrm{d}y,\quad \mbox{for all}\ u\in \mathcal{D}^{2,2}_0(\mathbb{R}^{N+2k}),
    \end{align}
    see \cite[(2.10)]{BWW03}, then combining with \eqref{rwevpbbbcb} and \eqref{rwevpbbbcbi}, we deduce that
    \begin{align*}
    & \left[(2^{**}-1)\Gamma_N-\Gamma_{N+2k}\right]
    \int_{\mathbb{R}^{N+2k}}
    (1+|y|^2)^{-4}|Z_{k,i}|^2 \mathrm{d}y
    \nonumber\\
    & \geq
    \frac{1}{2}(b^{-2}-1)\lambda_k
    [N(N-4)+(2\lambda_k-8)b^{-2}+2\lambda_k]
    \int_{\mathbb{R}^{N+2k}}\frac{|Z_{k,i}|^2}
    {|y|^{4}}\mathrm{d}y.
    \end{align*}
    Since
    \[
    (2^{**}-1)\Gamma_N\leq \Gamma_{N+2k},\quad \mbox{for all}\ k\geq 1,
    \]
    where $\Gamma_N$ is defined in \eqref{defgn}, then it holds that
    \begin{align}\label{rwevpbbbcbcf}
    (b^{-2}-1)\lambda_k
    [N(N-4)+(2\lambda_k-8)b^{-2}+2\lambda_k]
    \int_{\mathbb{R}^{N+2k}}\frac{|Z_{k,i}|^2}
    {|y|^{4}}\mathrm{d}y\leq 0.
    \end{align}
    Note that
    \[
    (b^{-2}-1)\lambda_k
    [N(N-4)+(2\lambda_k-8)b^{-2}+2\lambda_k]>0,\quad \mbox{for all}\ k\geq 1,
    \]
    due to $0<b<1$, $\lambda_k=k(N-2+k)$ and $N\geq 5$, then we conclude from \eqref{rwevpbbbcbcf} that $Z_{k,i}\equiv 0$ for all $k\geq 1$, that is, \eqref{rwevpb} does not exist nontrivial solutions for all $k\geq 1$, which means that we prove the claim. Note that this claim implies that $(\frac{N-4}{2}U_\mu+x\cdot\nabla U_\mu)$ is the only (up to multiplications) eigenfunction of $\nu_2= 2^{**}-1$.

    From the definition of eigenvalues, we deduce that $\nu_3>\nu_2$. Now, the proof of Theorem \ref{propep} is completed.
    \end{proof}

\vskip0.25cm

    \noindent{\bf \em Proof of Theorem \ref{thmlpnd}.}
    The proof directly follows from Theorem \ref{propep}.
    \qed

\vskip0.25cm

\section{{\bfseries Stability of Rellich-Sobolev inequality}}\label{sect:rt}

    In this section, we will prove the stability of Rellich-Sobolev inequality \eqref{RSi} and give the proof of Theorem \ref{thmprt}, inspired by Bianchi and Egnell \cite{BE91}.

    Since $\mathcal{M}_\mu=\{cU_{\mu,\lambda}: c\in\mathbb{R}, \lambda>0\}$ is two-dimensional manifold embedded in $\mathcal{E}_\mu$, and the tangential space at $(c,\lambda)$ is given by
    \begin{align*}
    T_{cU_{\mu,\lambda}}\mathcal{M}_\mu=\mathrm{Span}
    \left\{U_{\mu,\lambda},\ c\frac{\partial U_{\mu,\lambda}}{\partial \lambda}\right\},
    \end{align*}
    then Theorem \ref{propep} indicates that for any $u$ which is orthogonal to $T_{cU_{\mu,\lambda}}\mathcal{M}_\mu$,
    \begin{align}\label{czkj}
    \nu_3\int_{\mathbb{R}^N}U_{\mu,\lambda}^{2^{**}-2}u^2\mathrm{d}x \leq \|u\|^2_\mu,
    \end{align}
    where $\nu_3>2^{**}-1$ is independent of $\lambda$ given as in Theorem \ref{propep}.
    The main ingredient in the proof of Theorem \ref{thmprt} is contained in the lemma below, where the behavior near extremals set $\mathcal{M}_\mu$ will be studied.

    \begin{lemma}\label{lemma:rtnm2b}
    Suppose $N\geq 5$ and $0<\mu<N-4$. Then for any sequence $\{u_n\}\subset \mathcal{E}_\mu\backslash \mathcal{M}_\mu$ satisfying $\inf_n\|u_n\|_\mu>0$ and $\inf_{v\in \mathcal{M}_\mu}\|u_n-v\|_\mu\to 0$, it holds that
    \begin{equation}\label{rtnmb}
    \liminf\limits_{n\to\infty}\frac{\|u_n\|_\mu^2
    -\mathcal{S}_\mu\|u_n\|^2_{L^{2^{**}}(\mathbb{R}^N)}}
    {\inf\limits_{v\in \mathcal{M}_\mu}\|u_n-v\|^2_\mu}\geq 1-\frac{\nu_2}{\nu_3},
    \end{equation}
    where $\nu_2=2^{**}-1<\nu_3$ are given as in Theorem \ref{propep}.
    \end{lemma}

    \begin{proof}
    Let $d_n:=\inf_{v\in \mathcal{M}_\mu}\|u_n-v\|_\mu=\inf_{c\in\mathbb{R}, \lambda>0}\|u_n-cU_{\mu,\lambda}\|_\mu\to 0$. For each $n\in\mathbb{N}$ sufficiently large, there exist $c_n\in\mathbb{R}$ and $\lambda_n>0$ such that $d_n=\|u_n-c_nU_{\mu,\lambda_n}\|_\mu$. In fact,
    \begin{equation}\label{ikeda}
    \begin{split}
    \|u_n-cU_{\mu,\lambda}\|^2_\mu
    & =\|u_n\|^2_\mu+c^2\|U_{\mu,\lambda}\|^2_\mu-2c\langle u_n,U_{\mu,\lambda}\rangle_\mu \\
    & \geq \|u_n\|_\mu^2+c^2\|U_\mu\|^2_\mu-2|c|\|u_n\|_\mu \|U_\mu\|_\mu.
    \end{split}
    \end{equation}
    Thus, the minimizing sequence of $d_n^2$, say $\{c_{n,m},\lambda_{n,m}\}$, must satisfy $|c_{n,m}|\leq C$ for some $C\geq 1$ independent of $m$, which means that $\{c_{n,m}\}$ is bounded.
    On the other hand,
    \begin{align*}
    \left|\int_{|\lambda x|\leq \rho}\Delta u_n\Delta U_{\mu,\lambda} \mathrm{d}x\right|
    & \leq \int_{|y|\leq \rho}|\Delta (u_n)_{\frac{1}{\lambda}}(y)||\Delta U_\mu(y)| \mathrm{d}y \\
    & \leq \|\Delta u_n\|_{L^2(\mathbb{R}^N)}\left(\int_{|y|\leq \rho}|\Delta U_\mu|^2 \mathrm{d}y\right)^{1/2}
    = o_\rho(1)
    \end{align*}
    as $\rho\to 0$ which is uniformly for $\lambda>0$, where $(u_n)_{\frac{1}{\lambda}}(y)=\lambda^{-\frac{N-4}{2}}u_n(\lambda^{-1}y)$, and
    \begin{equation*}
    \left|\int_{|\lambda x|\geq \rho}\Delta u_n\Delta U_{\mu,\lambda} \mathrm{d}x \right|
    \leq \|\Delta U_\mu\|_{L^2(\mathbb{R}^N)}\left(\int_{|x|\geq \frac{\rho}{\lambda}}|\Delta u_n|^2 \mathrm{d}x\right)^{1/2}
    =  o_\lambda(1)
    \end{equation*}
    as $\lambda\to 0$ for any fixed $\rho>0$. By taking $\lambda\to 0$ and then $\rho\to 0$, we obtain
    \[
    \left|\int_{\mathbb{R}^N}\Delta u_n\Delta U_{\mu,\lambda} \mathrm{d}x\right| \to 0,\quad \mbox{as}\ \lambda\to 0.
    \]
    Similarly, we get
    \[
    \left|\int_{\mathbb{R}^N}\frac{\nabla u_n\cdot\nabla U_{\mu,\lambda}}{|x|^2} \mathrm{d}x\right| \to 0\quad \mbox{and}\quad \left|\int_{\mathbb{R}^N}\frac{u_n U_{\mu,\lambda}}{|x|^4} \mathrm{d}x\right| \to 0,\quad \mbox{as}\ \lambda\to 0.
    \]
    Therefore,
    \begin{small}
    \begin{align}\label{njt01}
    |\langle u_n,U_{\mu,\lambda}\rangle_\mu|
    & \leq \left|\int_{\mathbb{R}^N}\Delta u_n\Delta U_{\mu,\lambda} \mathrm{d}x\right|
    +C_{\mu,1}\left|\int_{\mathbb{R}^N}\frac{\nabla u_n\cdot\nabla U_{\mu,\lambda}}{|x|^2} \mathrm{d}x\right|
    +|C_{\mu,2}|\left|\int_{\mathbb{R}^N}\frac{u_n U_{\mu,\lambda}}{|x|^4} \mathrm{d}x\right|
    \nonumber \\
    & \to 0,\quad \mbox{as}\ \lambda\to 0.
    \end{align}
    \end{small}
    Moreover, by the explicit from of $U_{\mu,\lambda}$ we have
    \begin{equation*}
    \left|\int_{|\lambda x|\leq R}\Delta u_n\Delta U_{\mu,\lambda} \mathrm{d}x \right|
    \leq \|\Delta U_\mu\|_{L^2(\mathbb{R}^N)}\left(\int_{| x|\leq \frac{R}{\lambda}}|\Delta u_n|^2 \mathrm{d}x\right)^{1/2}
    = o_\lambda(1)
    \end{equation*}
    as $\lambda\to \infty$ for any fixed $R>0$, and
    \begin{equation*}
    \begin{split}
    \left|\int_{|\lambda x|\geq R}\Delta u_n\Delta U_{\mu,\lambda} \mathrm{d}x\right|
    & \leq \int_{|y|\geq R}|\Delta (u_n)_{\frac{1}{\lambda}}(y)||\Delta U_{\mu}(y)| \mathrm{d}y \\
    & \leq \|\Delta u_n\|_{L^2(\mathbb{R}^N)}\left(\int_{|y|\geq R}|\Delta U_\mu|^2 \mathrm{d}y\right)^{1/2}
    =  o_R(1)
    \end{split}
    \end{equation*}
    as $R\to \infty$ which is uniformly for $\lambda>0$. Thus, by taking $\lambda\to \infty$ and then $R\to \infty$, we also obtain
    \[
    \left|\int_{\mathbb{R}^N}\Delta u_n\Delta U_{\mu,\lambda} \mathrm{d}x\right| \to 0,\quad \mbox{as}\ \lambda\to \infty.
    \]
    Similarly, we deduce that
    \[
    \left|\int_{\mathbb{R}^N}\frac{\nabla u_n\cdot\nabla U_{\mu,\lambda}}{|x|^2} \mathrm{d}x\right| \to 0\quad \mbox{and}\quad \left|\int_{\mathbb{R}^N}\frac{u_n U_{\mu,\lambda}}{|x|^4} \mathrm{d}x\right| \to 0,\quad \mbox{as}\ \lambda\to \infty.
    \]
    Therefore,
    \begin{align}\label{njt02}
    |\langle u_n,U_{\mu,\lambda}\rangle_\mu|
    \to 0,\quad \mbox{as}\ \lambda\to \infty.
    \end{align}
    Combining with \eqref{njt01} and \eqref{njt02}, it follows from \eqref{ikeda} and $d_n\to 0$, $\inf_n\|u_n\|>0$ that the minimizing sequence $\{c_{n,m},\lambda_{n,m}\}$ must satisfy $1/C\leq |\lambda_{n,m}|\leq C$ for some $C\geq 1$ independent of $m$, which means that $\{\lambda_{n,m}\}$ is also bounded. Thus for each $n\in\mathbb{N}$ sufficiently large, $d_n$ can be attained by some $c_n\in\mathbb{R}$ and $\lambda_n>0$.

    Then we must have that $(u_n-c_n U_{\mu,\lambda_n})$ is perpendicular to $T_{c_n U_{\mu,\lambda_n}}\mathcal{M}_\mu$, i.e.,
    \[
    \langle U_{\mu,\lambda_n},u_n-c_n U_{\mu,\lambda_n}\rangle_\mu
    =\left\langle \frac{\partial U_{\mu,\lambda_n}}{\partial \lambda_n},u_n-c_n U_{\mu,\lambda_n}\right\rangle_\mu=0.
    \]
    Furthermore, the same as in \eqref{czkj}, we have
    \begin{equation}\label{epkeyibbg}
    \nu_3\int_{\mathbb{R}^N}U_{\mu,\lambda_n}^{2^{**}-2}(u_n-c_n U_{\mu,\lambda_n})^2\mathrm{d}x \leq \|u_n-c_n U_{\mu,\lambda_n}\|^2_\mu.
    \end{equation}
    Let $u_n=c_n U_{\mu,\lambda_n}+d_n w_n$, then $w_n$ is perpendicular to $T_{c_n U_{\mu,\lambda_n}}\mathcal{M}_\mu$,
    \[
    \|w_n\|_\mu=1\quad\mbox{and}\quad \|u_n\|^2_\mu=d_n^2+c_n^2\|U_\mu\|^2_\mu,
    \]
    in particular,
    \[
     \langle U_{\mu,\lambda_n},w_n\rangle_\mu
     =\int_{\mathbb{R}^N}U_{\mu,\lambda_n}^{2^{**}-1}w_n \mathrm{d}x
    =0.
    \]
    Then we can rewrite (\ref{epkeyibbg}) as follows:
    \begin{equation}\label{epkeyibbb}
    \int_{\mathbb{R}^N}U_{\mu,\lambda_n}^{2^{**}-2}w_n^2\mathrm{d}x
    \leq \frac{1}{\nu_3}.
    \end{equation}
    By  using Taylor's expansion, we deduce
    \begin{align}\label{epkeyiybb}
    \int_{\mathbb{R}^N}|u_n|^{2^{**}}\mathrm{d}x
    & = \int_{\mathbb{R}^N}|c_n U_{\mu,\lambda_n}+d_nw_n|^{2^{**}}\mathrm{d}x \nonumber\\
    & = |c_n|^{2^{**}}\int_{\mathbb{R}^N}U_{\mu,\lambda_n}^{2^{**}}\mathrm{d}x
    +d_n 2^{**}|c_n|^{2^{**}-1}
    \int_{\mathbb{R}^N}U_{\mu,\lambda_n}^{2^{**}-1}w_n \mathrm{d}x \nonumber\\
    & \quad +\frac{2^{**}(2^{**}-1)d_n^2  |c_n|^{2^{**}-2}}{2}
    \int_{\mathbb{R}^N}U_{\mu,\lambda_n}^{2^{**}-2}w_n^2\mathrm{d}x
    +o(d_n^2)  \nonumber\\
    & = |c_n|^{2^{**}}\|U_{\mu}\|_\mu^2
    + \frac{2^{**}(2^{**}-1)d_n^2  |c_n|^{2^{**}-2}}{2} \int_{\mathbb{R}^N}U_{\mu,\lambda_n}^{2^{**}-2}w_n^2\mathrm{d}x
    +o(d_n^2).
    \end{align}
    Then combining with (\ref{epkeyibbb}) and (\ref{epkeyiybb}), by the concavity of $t\mapsto t^{\frac{2}{2^{**}}}$, we obtain
    \begin{align}\label{epkeyiyxbb}
    \|u_n\|^2_{L^{2^{**}}(\mathbb{R}^N)}
    & \leq c_n^2\left[\|U_{\mu}\|^{2}_\mu+\frac{2^{**}(2^{**}-1)d_n^2  c_n^{-2}}{2\nu_3}
    +o(d_n^2)\right]^{\frac{2}{2^{**}}}  \nonumber\\
    & = c_n^2\left[\|U_{\mu}\|^{\frac{4}{2^{**}}}_\mu
    +\frac{2}{2^{**}}\cdot\frac{2^{**}(2^{**}-1)d_n^2  c_n^{-2}}{2\nu_3} \|U_{\mu}\|^{\frac{4}{2^{**}}-2}_\mu
    +o(d_n^2)\right] \nonumber\\
    & = c_n^2\|U_{\mu}\|^{\frac{4}{2^{**}}}_\mu+ \frac{d_n^2 (2^{**}-1)}{\nu_3}\|U_{\mu}\|^{\frac{4}{2^{**}}-2}_\mu+o(d_n^2).
    \end{align}
    Therefore, for $n$ sufficiently large,
    \begin{align*}
    & \|u_n\|_\mu^2
    -\mathcal{S}_\mu\|u_n\|^2_{L^{2^{**}}(\mathbb{R}^N)}
    \\
    & \geq d_n^2+c_n^2\|U_{\mu}\|^2_\mu
    - \mathcal{S}_\mu\left[c_n^2\|U_{\mu}\|^{\frac{4}{2^{**}}}_\mu+ \frac{d_n^2 (2^{**}-1)}{\nu_3}\|U_{\mu}\|^{\frac{4}{2^{**}}-2}_\mu+o(d_n^2)\right]  \nonumber\\
    & = d_n^2 \left(1-\frac{2^{**}-1}{\nu_3} \mathcal{S}_\mu \|U_{\mu}\|^{\frac{4}{2^{**}}-2}_\mu\right)
    +c_n^2\left(\|U_{\mu}\|^2_\mu- \mathcal{S}_\mu\|U_{\mu}\|^{\frac{4}{2^{**}}}_\mu\right)+o(d_n^2)  \\
    & = d_n^2\left(1-\frac{\nu_2}{\nu_3}\right)+o(d_n^2),
    \end{align*}
    using the fact that  $\|U_\mu\|_\mu^2
    =\mathcal{S}_\mu\left(\int_{\mathbb{R}^N}|U_\mu|^{2^{**}}
    \mathrm{d}x\right)^\frac{2}{2^{**}}$
     and $\|U_\mu\|_\mu^2
     =\int_{\mathbb{R}^N}|U_\mu|^{2^{**}}\mathrm{d}x$ which imply $\mathcal{S}_\mu= \|U_{\mu}\|^{2-\frac{4}{2^{**}}}_\mu$, then (\ref{rtnmb}) follows immediately.
    \end{proof}

    Now, we are ready to prove our main result.

\vskip0.25cm

    \noindent{\bf \em Proof of Theorem \ref{thmprt}.} We argue by contradiction. In fact, if the theorem is false then there exists a sequence $\{u_n\}\subset \mathcal{E}_\mu\setminus \mathcal{M}_\mu$ satisfying $\inf_n\|u_n\|_\mu>0$ and $\inf_{v\in \mathcal{M}_\mu}\|u_n-v\|^2_\mu\to 0$, such that
    \begin{equation*}
    \frac{\|u_n\|_\mu^2
    -\mathcal{S}_\mu\|u_n\|^2_{L^{2^{**}}(\mathbb{R}^N)}}
    {\inf\limits_{v\in \mathcal{M}_\mu}\|u_n-v\|^2_\mu}
    \to 0,\quad \mbox{as}\ n\to \infty.
    \end{equation*}
    By homogeneity, we can assume that $\|u_n\|_\mu=1$, and after selecting a subsequence we can assume that $\inf_{v\in \mathcal{M}_\mu}\|u_n-v\|_\mu\to \varpi\in[0,1]$ since $\inf_{v\in \mathcal{M}_\mu}\|u_n-v\|_\mu\leq \|u_n\|_\mu$. If $\varpi=0$, then we deduce a contradiction by Lemma \ref{lemma:rtnm2b}.

    The other possibility only is that $\varpi>0$, that is
    \[
    \inf_{v\in \mathcal{M}_\mu}\|u_n-v\|_\mu\to \varpi>0,\quad \mbox{as}\ n\to \infty,
    \]
    then we must have
    \begin{equation}\label{wbsi}
    \|u_n\|_\mu^2
    -\mathcal{S}_\mu\|u_n\|^2_{L^{2^{**}}(\mathbb{R}^N)}\to 0,\quad \|u_n\|_\mu=1.
    \end{equation}
    Since $\mathcal{S}_\mu<\mathcal{S}_0$ for all $0<\mu<N-4$, by taking the same arguments as those in \cite{WaWi00}, we can deduce that there is a sequence $\{\lambda_n\}\subset\mathbb{R}^+$ such that $\{(u_n)_{\lambda_n}\}$ contains a convergent subsequence which will lead to a contradiction.
    In fact, as stated in Remark \ref{remc}, we make the change of varable
    \begin{equation}\label{p2txyb}
    u_n(x)=|x|^{-\frac{\mu}{2}}v_n(|x|^{-\frac{\mu}{N-4}}x),
    \end{equation}
    then
    \begin{align*}
    0 \leq \|\Delta v_n\|^2_{L^2(\mathbb{R}^N)}
    -\mathcal{S}_0\|v_n\|^2_{L^{2^{**}}(\mathbb{R}^N)} \leq
    \left(1-\frac{\mu}{N-4}\right)^{-3}\left(\|u_n\|_\mu^2
    -\mathcal{S}_\mu\|u_n\|^2_{L^{2^{**}}(\mathbb{R}^N)}
    \right)
    \to 0.
    \end{align*}
    Note that \eqref{p2txyb} implies that $v_n$ does not invariant under translation, then by Lions' concentration-compactness principle (see \cite[Theorem \uppercase\expandafter{\romannumeral 2}.4]{Li85-1}), there is a sequence $\{\tau_n\}\subset\mathbb{R}^+$ such that
    \begin{equation*}
    \tau_n^{\frac{N-4}{2}}v_n(\tau_n x)\to W\quad \mbox{in}\ \mathcal{D}^{2,2}_0(\mathbb{R}^N),\quad \mbox{as}\ n\to \infty,
    \end{equation*}
    where $W(x)=c(d+|x|^2)^{-\frac{N-4}{2}}$ for some $c\neq 0$ and $d>0$. Then,
    \begin{equation*}
    \lambda_n^{\frac{N-4}{2}}u_n(\lambda_n x)\to U_*\quad \mbox{in}\ \mathcal{E}_\mu,\quad \mbox{as}\ n\to \infty,
    \end{equation*}
    for some $U_*\in\mathcal{M}_\mu$, where $\lambda_n=\tau_n^{\frac{N-4}{N-4-\mu}}$, which implies
    \begin{equation*}
    \inf\limits_{v\in \mathcal{M}_\mu}\|u_n-v\|_\mu
    =\inf\limits_{v\in \mathcal{M}_\mu}\left\|\lambda_n^{\frac{N-4}{2}}u_n(\lambda_n x)-v\right\|_\mu\to 0, \quad \mbox{as}\ n\to \infty,
    \end{equation*}
    this is a contradiction. Now, the proof of Theorem \ref{thmprt} is completed.
    \qed

\section{Weak-Lebesgue remainder term in a subdomain}\label{sectwlrt}

In this section, we will show a remainder term inequality in the weak Lebesgue-norm sense as in \eqref{thmprtbd}.
In order to prove \eqref{thmprtbd}, by homogeneity we can always assume that $\|u\|_{L^{2^{**}}(\Omega)}=1$. By H\"{o}lder's inequality we have
    \[
    \|u\|_{L^{\bar{p}}_w(\Omega)}\leq \|u\|_{L^{\bar{p}}(\Omega)}\leq |\Omega|^{\frac{2^{**}-\bar{p}}{2^{**}\cdot\bar{p}}}
    \|u\|_{L^{2^{**}}(\Omega)}
    =|\Omega|^{\frac{N-4-\mu}{2N}},
    \]
    where $\bar{p}=\frac{2N}{2(N-4)-\mu}$, thus it is sufficient to prove \eqref{thmprtbd} under $\|u\|_{\mu,\Omega}^2
    -\mathcal{S}_\mu
    \|u\|_{L^{2^{**}}(\Omega)}^2 \ll 1$. \eqref{thmprtbd} will follow immediately from the following lemma.

    \begin{lemma}\label{proprtbd}
    Suppose $N\geq 5$, $0<\mu<N-4$, and let $\Omega\subset \mathbb{R}^N$ with $|\Omega|<\infty$ which has smooth boundary and contains origin in its interior. Then there exists a constant $\mathcal{B}=\mathcal{B}(N,\mu)>0$ independent of $\Omega$ such that for all $u\in \dot{\mathcal{E}}_{\mu,\Omega}$ satisfying $\|u\|_{L^{2^{**}}(\Omega)}=1$ and $\|u\|_{\mu,\Omega}^2
    -\mathcal{S}_\mu
    \|u\|_{L^{2^{**}}(\Omega)}^2 \ll 1$,
    \begin{equation}\label{rtbdke}
    \|u\|_{L^{\bar{p}}_w(\Omega)}
    \leq \mathcal{B} |\Omega|^{\frac{N-4-\mu}{N}}\inf_{v\in \mathcal{M}_\mu}\|u-v\|_{\mu},
    \end{equation}
    where $L^{\bar{p}}_w(\Omega)$ denotes the weak $L^{\bar{p}}$-norm given as in \eqref{defwln}.
    \end{lemma}

    \begin{proof}
    We follow the arguments as those in \cite[Proposition 3]{CFW13}. Firstly, Theorem \ref{thmprt} and $\|u\|_{\mu,\Omega}^2
    -\mathcal{S}_\mu
    \|u\|_{L^{2^{**}}(\Omega)}^2 \ll 1$ indicate
    \[
    \inf_{v\in \mathcal{M}_\mu}\|u-v\|_{\mu}\ll 1,
    \]
    then as in the proof of Lemma \ref{lemma:rtnm2b}, we have
    \[
    \inf_{v\in \mathcal{M}_\mu}\|u-v\|_{\mu}=\|u-cU_{\mu,\lambda}\|_{\mu},\quad  \mbox{for some}\ c\in\mathbb{R}, \ \lambda>0.
    \]
    Note that $\|u\|_{L^{2^{**}}(\Omega)}=1$ and $\|u\|_{\mu,\Omega}^2
    -\mathcal{S}_\mu
    \|u\|_{L^{2^{**}}(\Omega)}^2 \ll 1$ imply that $\|u\|_{\mu,\Omega}$ is bounded away from zero and infinity, i.e., $c_0\leq \|u\|_{\mu,\Omega}\leq C_0$ for some constants $C_0\geq c_0>0$.  In fact, we can always deduce
    \[
    \mathcal{S}^{1/2}_\mu\leq \|u\|_{\mu,\Omega}\leq 2\mathcal{S}^{1/2}_\mu.
    \]

    Let $\tau=\tau(N,\mu)\in (0,\mathcal{S}^{1/2}_\mu)$ be given by
    \begin{align}\label{defrho}
    \frac{\tau\|U_{\mu}\|_{L^{2^{**}}(\mathbb{R}^N)}}
    {\mathcal{S}^{1/2}_\mu-\tau}
    =\|U_{\mu}\|_{L^{2^{**}}(\mathbb{R}^N\setminus \mathrm{B}_1)}
    =K_{N,\mu}\left(|\mathbb{S}^{N-1}|\int^\infty_1 \frac{r^{N-1-\frac{N\mu}{N-4}}}
    {\left[1+r^{2(1-\frac{\mu}{N-4})}\right]^N}\mathrm{d}r
    \right)^{1/2^{**}},
    \end{align}
    where $\mathrm{B}_1:=\mathrm{B}(\mathbf{0},1)$. It is obvious that $\tau$ is a fixed constant. So
    \[
    \inf_{v\in \mathcal{M}_\mu}\|u-v\|_{\mu}\leq \tau,
    \]
    thanks to $\inf_{v\in \mathcal{M}_\mu}\|u-v\|_{\mu}\ll 1$. Note that
    \begin{align*}
    |\|u\|_{\mu}-\|cU_{\mu,\lambda}\|_{\mu}|
    \leq \|u-cU_{\mu,\lambda}\|_{\mu}
    =\inf_{v\in \mathcal{M}_\mu}\|u-v\|_{\mu}\leq \tau,
    \end{align*}
    which implies that
    \[
    \frac{\mathcal{S}^{1/2}_\mu-\tau}{\|U_{\mu}\|_{\mu}}
    \leq \frac{\|u\|_{\mu}-\tau}{\|U_{\mu}\|_{\mu}}\leq |c|\leq \frac{\|u\|_{\mu}+\tau}{\|U_{\mu}\|_{\mu}}
    \leq \frac{2\mathcal{S}^{1/2}_\mu+\tau}{\|U_{\mu}\|_{\mu}}.
    \]
    Then we have
    \begin{align}\label{dug}
    \inf_{v\in \mathcal{M}_\mu}\|u-v\|_{\mu}^2
    & =\|u-cU_{\mu,\lambda}\|_{\mu}^2
    \geq \mathcal{S}_\mu \|u-cU_{\mu,\lambda}\|^2_{L^{2^{**}}(\mathbb{R}^N)}
    \geq \mathcal{S}_\mu |c|^2\|U_{\mu,\lambda}\|^2_{L^{2^{**}}(\mathbb{R}^N\setminus\Omega)}
    \nonumber\\
    & \geq \left(\frac{\mathcal{S}^{1/2}_\mu-\tau}{\|U_{\mu}\|_{\mu}}\right)^2
    \mathcal{S}_\mu \|U_{\mu,\lambda}\|^2_{L^{2^{**}}(\mathbb{R}^N\setminus\Omega)}
    = \frac{(\mathcal{S}^{1/2}_\mu-\tau)^2}
    {\|U_{\mu}\|^2_{L^{2^{**}}(\mathbb{R}^N)}} \|U_{\mu,\lambda}\|^2_{L^{2^{**}}(\mathbb{R}^N\setminus\Omega)}.
    \end{align}
    Now let $\mathrm{B}_R\subset\mathbb{R}^N$ denote the open ball at origin with radius $R$ such that $|\mathrm{B}_R|=|\Omega|$. Since $U_{\mu,\lambda}$ is radial and strictly decreasing, the bathtub principle \cite[Theorem 1.14]{LL01} implies that
    \begin{align*}
    \|U_{\mu,\lambda}\|^2_{L^{2^{**}}(\mathbb{R}^N\setminus\Omega)}
    \geq \|U_{\mu,\lambda}\|^2_{L^{2^{**}}(\mathbb{R}^N\setminus \mathrm{B}_R)},
    \end{align*}
    then
    \begin{align}\label{Uwgo}
    \|U_{\mu,\lambda}\|_{L^{2^{**}}(\mathbb{R}^N\setminus \mathrm{B}_R)}^{2^{**}}
    & \leq \left(\frac{\|U_{\mu}\|_{L^{2^{**}}(\mathbb{R}^N)}\inf\limits_{v\in \mathcal{M}_\mu}\|u-v\|_{\mu}}
    {\mathcal{S}^{1/2}_\mu-\tau}\right)^{2^{**}}
    \nonumber\\
    & \leq \left(\frac{\tau\|U_{\mu}\|_{L^{2^{**}}(\mathbb{R}^N)}}
    {\mathcal{S}^{1/2}_\mu-\tau}\right)^{2^{**}}
    =K_{N,\mu}^{2^{**}}|\mathbb{S}^{N-1}|\int^\infty_1 \frac{r^{N-1-\frac{N\mu}{N-4}}}
    {\left[1+r^{2(1-\frac{\mu}{N-4})}\right]^N}\mathrm{d}r
    \end{align}
    by our choice of $\tau$ in \eqref{defrho}. On the other hand, we compute
    \begin{align*}
    \|U_{\mu,\lambda}\|_{L^{2^{**}}(\mathbb{R}^N\setminus \mathrm{B}_R)}^{2^{**}}
    & =K_{N,\mu}^{2^{**}}|\mathbb{S}^{N-1}|\int^\infty_R \frac{\lambda^{\frac{(N-4-\mu)N}{N-4}}r^{N-1-\frac{N\mu}{N-4}}}
    {\left[1+(\lambda r)^{2(1-\frac{\mu}{N-4})}\right]
    ^N}\mathrm{d}r
    \\
    & =K_{N,\mu}^{2^{**}}|\mathbb{S}^{N-1}|\int^\infty_{\lambda R} \frac{r^{N-1-\frac{N\mu}{N-4}}}
    {\left[1+r^{2(1-\frac{\mu}{N-4})}\right]^N}\mathrm{d}r,
    \end{align*}
    which implies that $\lambda R\geq 1$, therefore
    \begin{align}\label{Uwg}
    \|U_{\mu,\lambda}\|_{L^{2^{**}}(\mathbb{R}^N\setminus \mathrm{B}_R)}^{2^{**}}
    & \geq 2^{-N}K_{N,\mu}^{2^{**}}|\mathbb{S}^{N-1}|
    \int^\infty_{\lambda R} r^{-N-1+\frac{N\mu}{N-4}}\mathrm{d}r
    \nonumber\\
    & =2^{-N}K_{N,\mu}^{2^{**}}|\mathbb{S}^{N-1}|
    \frac{N-4}{(N-4-\mu)N}
    R^{-\frac{(N-4-\mu)N}{N-4}}
    \lambda^{-\frac{(N-4-\mu)N}{N-4}}.
    \end{align}
    Combining \eqref{Uwgo} and \eqref{Uwg}, we conclude from \eqref{dug} that
    \begin{align}\label{dugf}
    \inf_{v\in \mathcal{M}_\mu}\|u-v\|_{\mu}\geq C_1R^{-\frac{N-4-\mu}{2}}\lambda^{-\frac{N-4-\mu}{2}}
    \end{align}
    with
    \[
    C_1:=\frac{\mathcal{S}^{1/2}_\mu-\tau}
    {\|U_{\mu}\|_{L^{2^{**}}(\mathbb{R}^N)}}
    K_{N,\mu}\left(2^{-N}|\mathbb{S}^{N-1}|
    \frac{N-4}{(N-4-\mu)N}\right)^{1/2^{**}}.
    \]
    Then from \eqref{defrho} and \eqref{dugf}, we have
    \begin{align}\label{rtbdlbcdl}
    \|u\|_{L^{\bar{p}}_w(\Omega)}
    & \leq \|cU_{\mu,\lambda}\|_{L^{\bar{p}}_w(\Omega)}
    + \|u-cU_{\mu,\lambda}\|_{L^{\bar{p}}_w(\Omega)}
    \nonumber \\
    & \leq |c|\lambda^{-\frac{N-4-\mu}{2}}
    \|U_\mu\|_{L^{\bar{p}}_w(\mathbb{R}^N)}
    + \|u-cU_{\mu,\lambda}\|_{L^{\bar{p}}(\Omega)}
    \nonumber \\
    & \leq \frac{2\mathcal{S}^{1/2}_\mu+\tau}
    {\|U_{\mu}\|_{\mu}}
    \lambda^{-\frac{N-4-\mu}{2}}
    \|U_\mu\|_{L^{\bar{p}}_w(\mathbb{R}^N)}
    +|\Omega|^{\frac{N-4-\mu}{2N}}\mathcal{S}_\mu^{-1/2}
    \|u-cU_{\mu,\lambda}\|_{\mu}
    \nonumber \\
    & \leq \mathcal{B}|\Omega|^{\frac{N-4-\mu}{2N}}
    \inf_{v\in \mathcal{M}_\mu}\|u-v\|_{\mu}
    \end{align}
    with
    \[
    \mathcal{B}
    :=\frac{(2\mathcal{S}^{1/2}_\mu+\tau)
    \|U_\mu\|_{L^{\bar{p}}_w(\mathbb{R}^N)}}
    {C_1|\mathbb{S}^{N-1}|^{\frac{2N}{N-4-\mu}}
    \|U_{\mu}\|_{\mu}}
    +\mathcal{S}_\mu^{-1/2},
    \]
    where we have used the fact that $|\Omega|=|\mathrm{B}_R|=|\mathbb{S}^{N-1}| R^N$. It is not difficult to verify that $\|U_\mu\|_{L^{\bar{p}}_w(\mathbb{R}^N)}<\infty$.
    Now, the proof of \eqref{rtbdke} is completed.
    \end{proof}

Now, we are ready to give the proof of weak-Lebesgue remainder term inequality \eqref{thmprtbd}.
    \vskip0.25cm

\noindent{\bf \em Proof of Theorem \ref{thmrtbd}.} The proof directly follows from Theorem \ref{thmprt} and Lemma \ref{proprtbd}.
\qed

\vskip0.25cm

\noindent{\bfseries Data availability}
No data was used for the research described in the article.

\vskip0.25cm

\noindent{\bfseries Acknowledgements} We are very grateful to the anonymous referees for insightful suggestions, which lead to an important improvement of the paper.
The research has been supported by National Natural Science Foundation of China (No. 12371121).

    \end{document}